\documentclass[12pt]{article}
\usepackage[dvips]{graphicx}
\usepackage{amsmath}
\usepackage{amsfonts}
\usepackage{amssymb}
\newtheorem{theorem}{Theorem}[subsection]

\newtheorem{corollary}[theorem]{Corollary}

\newtheorem{definition}[theorem]{Definition}
\newtheorem{example}[theorem]{Example}

\newtheorem{lemma}[theorem]{Lemma}

\newtheorem{problem}[theorem]{Problem}
\newtheorem{proposition}[theorem]{Proposition}

\newenvironment{proof}[1][Proof]{\textbf{#1.} }{\ \rule{0.5em}{0.5em}}

\begin{document}

\title{Signed permutations and the four color theorem}
\author{Shalom Eliahou and C\'{e}dric Lecouvey\\Laboratoire de Math\'{e}matiques Pures et Appliqu\'{e}es \\Universit\'{e} du Littoral C\^{o}te d'Opale\\50 rue F. Buisson, B.P. 699\\62228 Calais Cedex, France}
\date{}
\maketitle
\begin{abstract}
To each permutation $\sigma$ in $S_{n}$ we associate a triangulation of a
fixed $(n+2)$-gon.\ We then determine the fibers of this association and show
that they coincide with the sylvester classes depicted in \cite{HNT}. A signed
version of this construction allows us to reformulate the four color theorem
in terms of the existence of a signable path between any two permutations in
the Cayley graph of the symmetric group $S_{n}.$
\end{abstract}

\baselineskip5.5mm

\section{Introduction}

In this paper, we obtain a reformulation of the four color theorem in terms of
signed permutations (Theorem \ref{fc}).\ Signed permutations are standard
words (that is with no letter repeated) $w$ of length $n$ on the alphabet
$\{\overline{n},\ldots,\overline{1},1,\ldots,n\}$ which do not contain any
pair $(k,\overline{k}).$ The barred (resp.\ unbarred) letters are interpreted
as negative (resp. positive) integers and we write $\left|  w\right|  $ for
the word on $\{1,\ldots,n\}$ obtained by erasing the bars which appear on the
letters of $w.\;$We first describe a combinatorial procedure which attaches to
each permutation $\sigma\in S_{n}$ (considered as a standard word on a given
totally ordered alphabet $X_{n}$) a triangulation $\varphi(\sigma)$ of a fixed
$(n+2)$-gon $P$.\ The fibers of this association coincide with the sylvester
classes defined in \cite{HNT}. Our construction is close to that used by
Reading in \cite{Re} and the sylvester congruence on \textit{standard} words
is in fact a special case of cambrian congruence.\ Nevertheless, we have
chosen to give an independent exposition of the occurrence of the sylvester
relations in the context of triangulations.\ This makes the paper
self-contained and permits us to expose precisely the results which are
required in our reformulation of the four color theorem in terms of signed
permutations. This connection with the results of \cite{Re} motivates an
additional reformulation of the four color theorem in terms of the geometry of
the associahedron.

We also introduce colored triangulations obtained by associating a color
(considered as a letter of a totally ordered alphabet $\mathcal{C}$ with $p$
colors) to each face (or to each vertex) of a given triangulation of $P$.\ The
map $\varphi$ admits a natural extension $\Phi$ defined from the set of words
$w$ with letters in $\mathcal{C}$ to a particular subset of colored
triangulations we have called simple.\ We show that the fiber $\Phi^{-1}\{w\}$
coincides with the sylvester class of $w$ (Theorem \ref{Th_sylvT}).\ This
implies that the simple colored triangulations can be regarded as
combinatorial objects analogous to binary search trees.\ The map $\Phi$ can be
interpreted as an insertion scheme for simple triangulations.\ This gives an
alternative to Knuth's insertion algorithm on binary search trees.

\textit{Signed triangulations} are colored triangulations with $\mathcal{C}%
=\{\pm1\}$.\ To each triangulation $T$ of $P$ corresponds the sylvester class
$\varphi^{-1}(T).\;$This permits us to associate to each signed triangulation
$T_{\varepsilon},$ defined as a signing $\varepsilon$ of $T$, the set
$[T]_{s}$ of signed permutation words $w$ such that $\left|  w\right|
\in\varphi^{-1}(T)$.\ To obtain our reformulation of the four color theorem,
we establish that for each signed flip (defined as a particular diagonal flip)
in $T_{\varepsilon}$ one can find two words $w_{1}$ and $w_{2}$ in $[T]_{s}$
such that%
\[%
\begin{array}
[c]{rcl}%
w_{1} & = & u\ \alpha\ \gamma\ v\\
w_{2} & = & u\ \overline{\gamma}\ \overline{\alpha}\ v
\end{array}
\]
where $\alpha,\gamma$ are letters with the same sign, $u,v$ are factors of
$w_{1}$ and $v$ does not contain any letter $\beta$ such that $\left|
\alpha\right|  <\left|  \beta\right|  <\left|  \gamma\right|  .$

We also study the combinatorial problem of describing the graph obtained from
a colored triangulation $T_{\varepsilon}$ by computing successive flips.\ When
no condition is imposed on the flip, it is well known that the flip graph of
the (ordinary) triangulations $T$ is connected and contains all the
triangulations of $P.\;$In addition to signed flips, we also consider in this
paper switched flips.\ They are flips which preserve the coloring and for
which a flip operation is authorized when the two faces considered have
distinct colors.\ For any $\mu\in\mathbb{N}^{p}$ such that $\mu_{1}+\cdot
\cdot\cdot+\mu_{p}=n$, let $S_{\mu}$ be the Frobenius subgroup of $S_{n}$
defined by $\mu.\;$We prove that the map $\Phi$ is a morphism from the Cayley
graph of $S_{n}/S_{\mu}$ to the graph whose vertices are the simple colored
triangulations associated to $\mu$ connected by an edge when they differ by a
switched flip. So the graph defined from a simple colored triangulation by
applying switched flips contains all the simple colored triangulations which
have the same coloring and thus, is connected.

\bigskip

The paper is organized as follows. Section $2$ is devoted to the combinatorial
background on triangulations, colored triangulations and flip operations we
need in the sequel. In section $3,$ we introduce the maps $\varphi$ and $\Phi$
and prove that their fibers coincide with sylvester classes.\ The
reformulations of the four color theorem using respectively the signed
permutations and the geometry of the associahedron are given in Section
$4$.\ Finally, we study in section $5$ the flip graphs generated by colored
triangulations when restrictive conditions are imposed on the flips.

\section{Triangulations and flips\label{subsec_label}}

\subsection{Triangulations of an $(n+2)$-gon}

Consider an integer $n\geq0$ and $X=\{x_{1}<\ldots< x_{n}\}$ a subset of
$\mathbb{N} \setminus\{0\}$ of size $n.\;$We denote by $\widehat{X}$ the
augmented set $\widehat{X} = X \cup\{0, \infty\}$, totally ordered as
follows:
\[
\widehat{X}=\{0<x_{1}<\ldots<x_{n}<\infty\}.
\]
Each integer $x\in X$ has a predecessor and a successor in $\widehat{X}$ that
we denote $\mathrm{pred}_{\widehat{X}}(x)$ and $\mathrm{succ}_{\widehat{X}%
}(x)$ (or $\mathrm{pred}(x)$ and $\mathrm{succ}(x)$, for short), respectively.
Thus,
\[
\mathrm{pred}(x) < x < \mathrm{succ}(x)
\]
is a 3-element interval in $\widehat{X}$. \ Let $P=P_{n}$ be a convex
$(n+2)$-gon, with vertices labelled by the elements of $\widehat{X}$ in a
clockwise increasing way. By a \textit{triangulation} $T$ of $P$, we mean a
plane graph with the $n+2$ vertices and $n+2$ edges of $P$, and with $n-1$
additional edges, called \textit{diagonals}, which subdivide the inner face of
$P$ into $n$ triangles, called the \textit{faces} of $T$.

We denote by $\mathcal{T}_{n}$ the set of triangulations of the polygon
$P=P_{n}.\;$It is well known that the cardinality of $\mathcal{T}_{n}$ is
equal to the Catalan number $c_{n}=\frac{1}{2n+1}\binom{2n}{n}.$ By a slight
abuse of notation, we shall make no distinction between a vertex $v$ of $P$
and its label in $\widehat{X}$. In other words, we consider the polygon $P$
and its triangulations $T$ as graphs with vertex set $\widehat{X}.$

Recall that the \textit{degree} of a vertex $v$ in a graph $G$ is the number
of edges of $G$ which are incident to $v.\;$ Now, in any triangulation $T$ of
the polygon $P$, as above, the degree of each vertex $v$ is at least 2,
because of the two edges of $P$ which are incident to $v$. An \textit{ear} in
$T$ is a vertex of degree exactly $2$, i.e. a vertex incident to zero
diagonals of $T$. Thus, an ear $e$ belongs to exactly one face of $T$, whose
three vertices are $\{\mathrm{pred}(e),e,\mathrm{succ}(e)\}.$ (More generally,
a vertex of degree $d$ in $T$ belongs to exactly $d-1$ faces of $T$.) It is
easy to see that every triangulation of $P$ contains at least two ears, and
that no two ears are adjacent if $n\geq2$.

By convention, we shall label each face $F$ of a triangulation $T$ by its
middle vertex. In other words, if $x < y < z \in\widehat{X}$ are the three
vertices of the face $F$, we shall label $F$ by $y$. This labelling gives a
bijection between the faces of $T$ and the set $X$; there are no faces
labelled $0$ of $\infty$.%

\begin{center}
\includegraphics[
height=4.2878cm,
width=4.0857cm
]%
{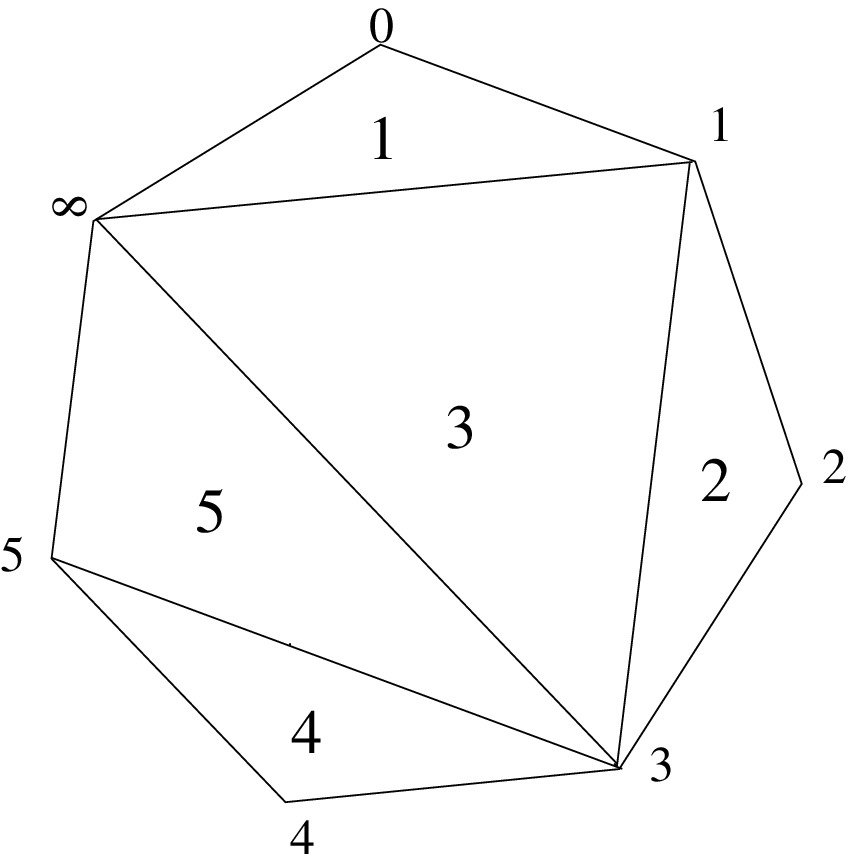}%
\\
Figure 1: Labelling of the faces
\label{figure2}%
\end{center}

\subsection{Colored triangulations}

Consider now a totally ordered set $\mathcal{C}=\{c_{1}<\ldots<c_{p}\}$, which
will be referred to in the sequel as the \textit{set of colors}.\ Given a
triangulation $T$ and $\varepsilon=(\varepsilon_{1},\ldots,\varepsilon_{n})
\in\mathcal{C}^{n}$ a $n$-tuple of colors, we call \textit{colored
triangulation} the triangulation $T_{\varepsilon}$ obtained by replacing in
$T$ each label $x_{i}\in X$ by the pair $(x_{i},\varepsilon_{i}).$ The faces
of $T_{\varepsilon}$ are colored by associating the color $\varepsilon_{i}$ to
the face labelled $x_{i}$ by the previous procedure. The vertex $x_{i}$ and
the face $T_{i}$ in the triangulation $T_{\varepsilon}$ are hence colored by
the color $\varepsilon_{i}\in\mathcal{C}$. Note that the additional vertices
labelled by $0$ and $\infty$ are not colored. We write simply $T$ for the
underlying triangulation associated to $T_{\varepsilon}$.\ Denote by
$\mathcal{T}_{n}(\varepsilon)$ the set of triangulations of $P$ colored by
$\varepsilon.\;$When $p=2,$ we set $\mathcal{C}=\{-,+\}$ (with $- < +$) and
call \textit{signed triangulations} the colored triangulations on
$\mathcal{C}.$ We will say that $\varepsilon$ is an \textit{increasing
coloring} when $\varepsilon_{1}\leq\cdot\cdot\cdot\leq\varepsilon_{n}.$

\begin{definition}
\label{de_coltrian}The colored triangulation $T_{\varepsilon}$ is called
simple if it verifies the following conditions:

\begin{itemize}
\item $\varepsilon$ is an increasing coloring;

\item  there are no inner diagonals in $T_{\varepsilon}$ connecting two
vertices with the same color;

\item  when two consecutive vertices $x_{i}$ and $x_{i+1}$ are colored with
the same color, the third vertex $t_{i}$ of the face they define verifies
$t_{i}<x_{i}.$
\end{itemize}
\end{definition}

\noindent\textbf{Remarks.}

\noindent$\mathrm{(i)}$ In a simple colored triangulation, there is no face
with three vertices of the same color. Similarly, if a face has two vertices
colored with the same color, these vertices must be consecutive.

\noindent$\mathrm{(ii)}$ In a colored triangulation, the coloring of the
vertices determines that of the faces.\ Conversely, the coloring of the faces
gives the coloring of the vertices since each face of a triangulation has the
color of its middle vertex. So the colored triangulation $T_{\varepsilon}$ is
characterized by $T$ and one of these two colorings.%

\begin{gather}%
\raisebox{-0cm}{\includegraphics[
height=4.7886cm,
width=5.0918cm
]%
{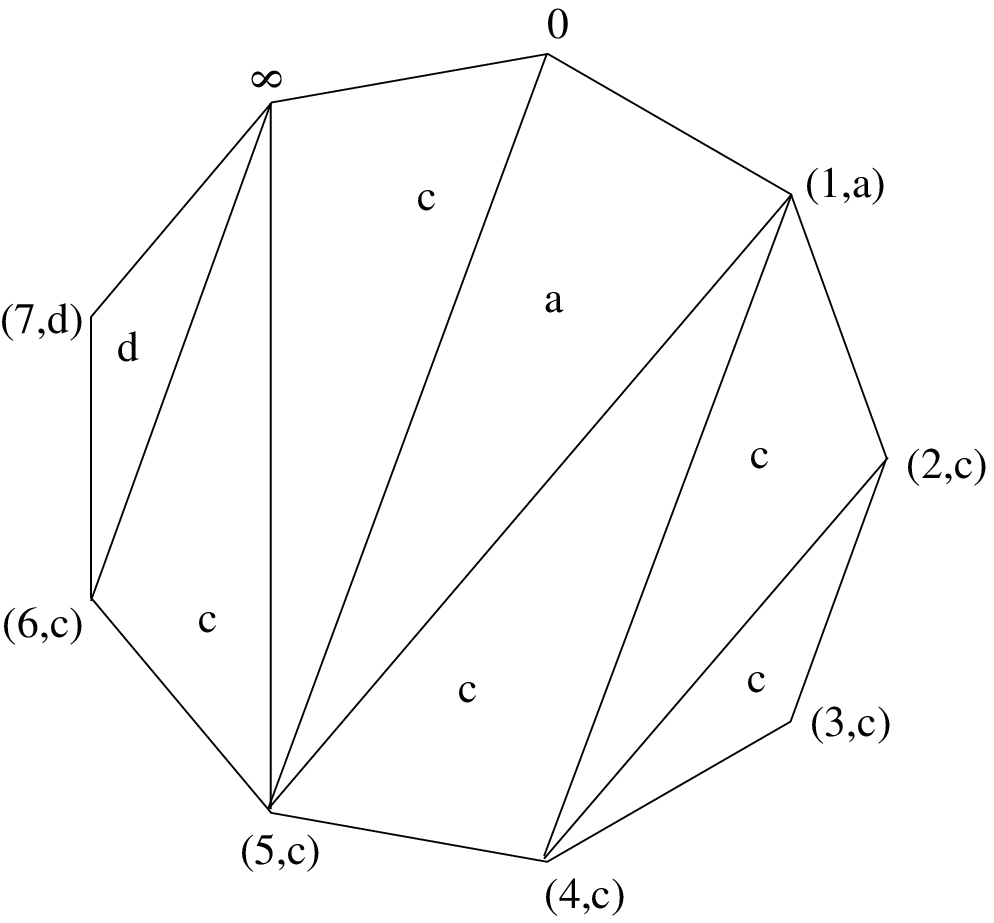}%
}%
\raisebox{-0cm}{\includegraphics[
height=4.8084cm,
width=5.0918cm
]%
{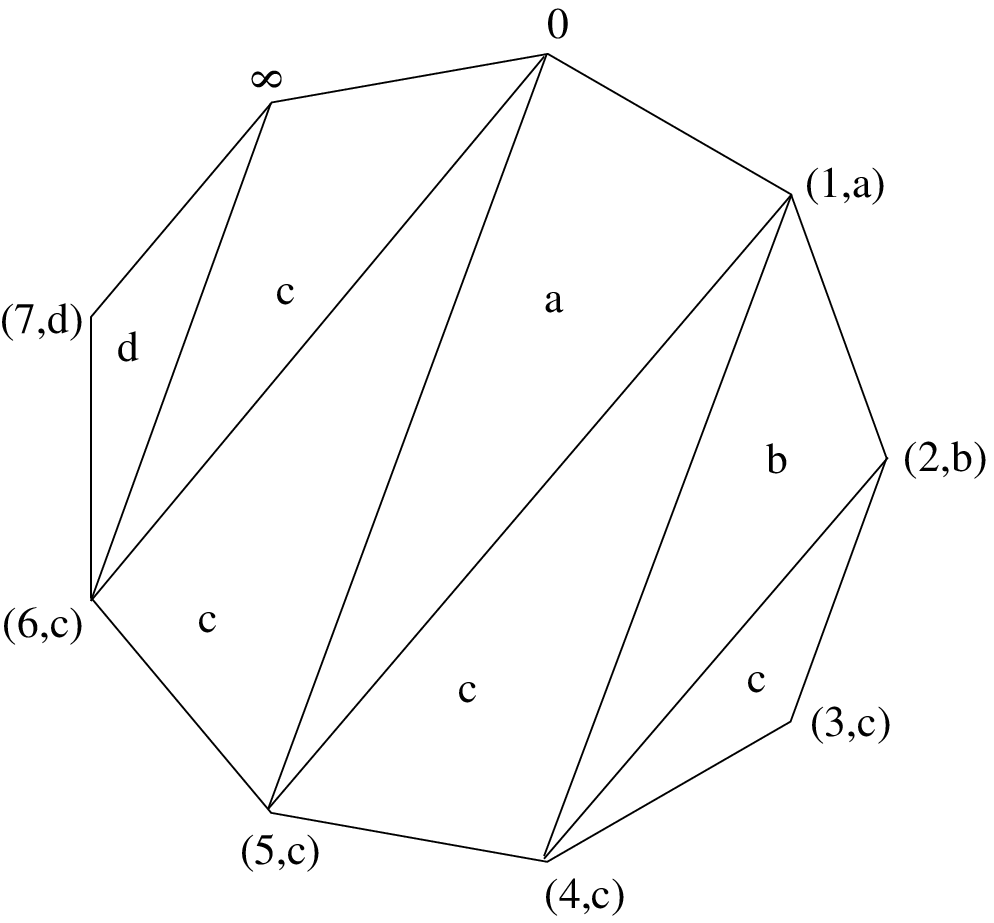}%
}%
\label{vc}\\
\text{Figure 2: Non simple and simple colored triangulations for }%
\mathcal{C}=\{a<b<c\}\nonumber
\end{gather}
The set of simple triangulations colored by $\varepsilon$ will be denoted by
$\mathcal{ST}_{n}(\varepsilon).$

\subsection{Flips and restricted flips on colored triangulations}

Consider $T\in\mathcal{T}_{n}$ a triangulation of the polygon $P=P_{n}.$ The
triangulation $T$ may be transformed into another one by a diagonal flip. The
\textit{diagonal flip}, or \textit{flip} for short, of a diagonal $d$ in $T$
is the following operation: in the quadrilateral $Q=v_{1}v_{2}v_{3}v_{4}$
formed by the two faces of $T$ adjacent to $d=v_{2}v_{4}$, remove $d$ and
replace it by the opposite diagonal $d^{\prime}=v_{1}v_{3}.$ The result is a
new triangulation $T^{\prime}$ of $P$ (see figure below). It is well known
that given any triangulations $T_{1}$ and $T_{2}$ of $P,$ there exists a
sequence of diagonal flips transforming $T_{1}$ into $T_{2}.\;$The graph on
$\mathcal{T}_{n}$ obtained by joining two triangulations which differ by
exactly one flip is called the \textit{flip graph} (see \cite{BW}).\ We will
denote it by $\mathcal{F}_{n}.\;$There are many other labellings of the flip
graph by objects enumerated by the Catalan numbers, such as binary trees,
parenthesizations, etc.

\bigskip

We now introduce various types of \textit{colored flips} on colored
triangulations of $P$, by imposing special constraints on the two faces
adjacent to the diagonal being flipped. The notion of signed flips plays a key
role in our reformulation of the four color theorem.

\begin{itemize}
\item  A \textit{signed flip} in a signed triangulation is a diagonal flip
such that the two faces adjacent to the flipped diagonal have equal signs,
which are changed after the flip.

\item  A \textit{homogeneous flip} in a colored triangulation is a diagonal
flip such that the two faces adjacent to the flipped diagonal have the same
color, which are preserved after the flip.

\item  A \textit{switched flip} is a diagonal flip which preserves the vertex
colors and such that the two faces adjacent to the flipped diagonal have
different colors.
\end{itemize}

Such flips will be referred to in the sequel as \textit{restricted flips}.%

\begin{center}
\includegraphics[
height=10.0276cm,
width=12.1166cm
]%
{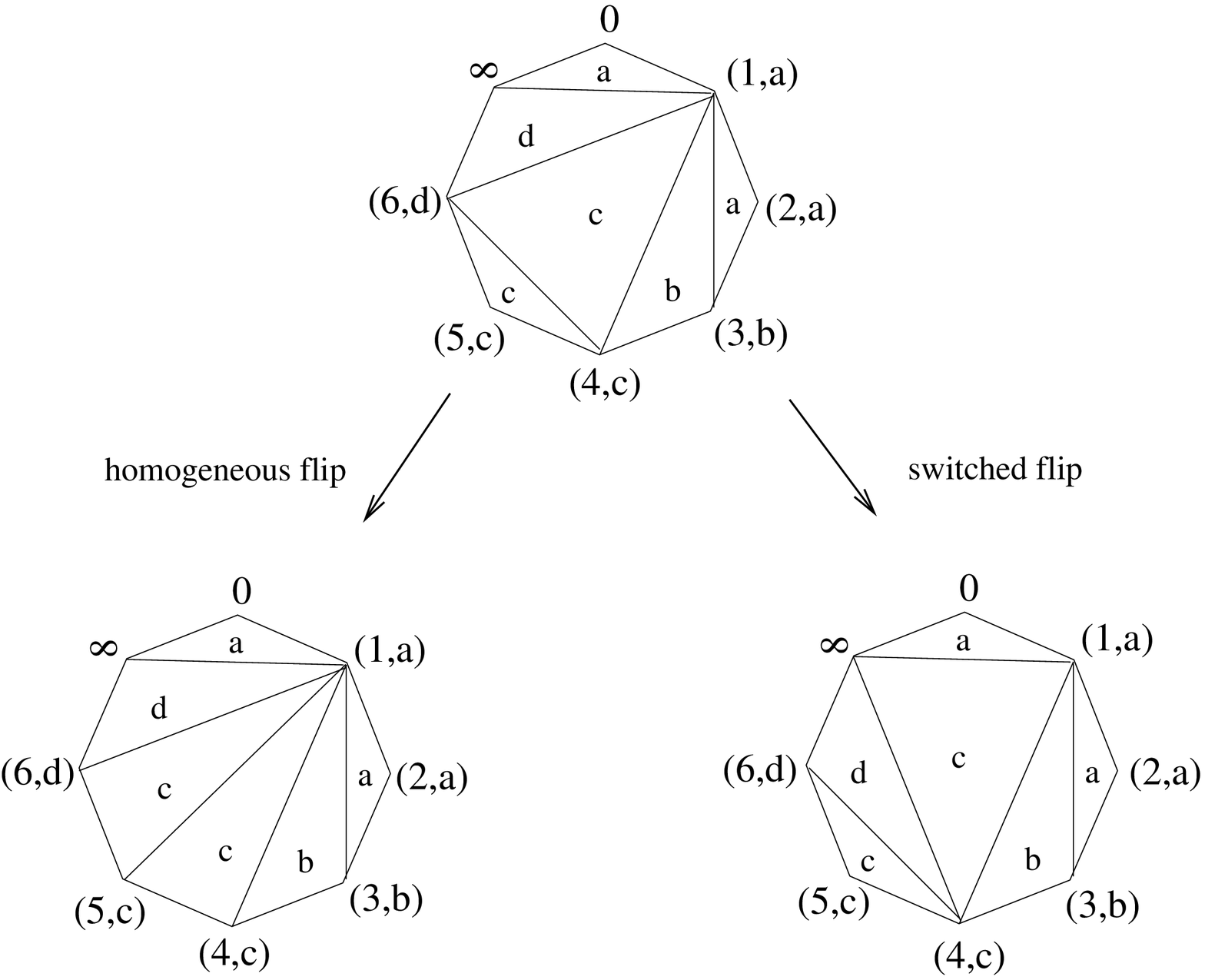}%
\label{fig6}%
\end{center}

\begin{center}
\includegraphics[
height=4.5624cm,
width=12.1166cm
]%
{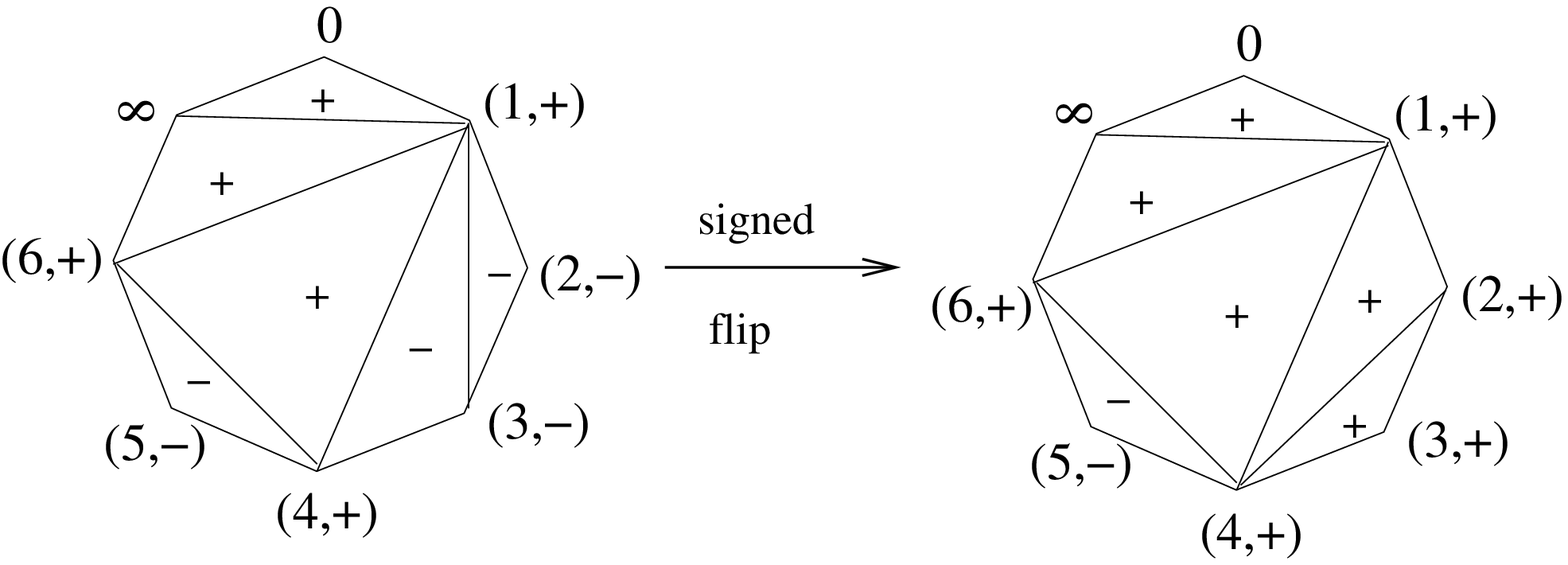}%
\\
Figure 3
\label{fig3}%
\end{center}

%

\section{Triangulations and sylvester relations\label{sec_sylv}}

\subsection{The triangulation associated to a permutation}

As above, let $X=\{x_{1}<x_{2}<\ldots<x_{n}\}$ be a subset of $\mathbb{N}%
\setminus\{0\}$. As usual, we denote by $S_{n}$ the symmetric group of rank
$n$.\ By realizing $S_{n}$ as the permutation group of $X,$ one can identify
each $\sigma\in S_{n}$ with a standard word $\sigma=x_{i_{1}}\ldots x_{i_{n}}$
on $X.$ We shall now define a map
\[
\varphi:S_{n}\rightarrow\mathcal{T}_{n}%
\]
which, being surjective, will allow us to represent each triangulation of $P$
by a suitable word in $S_{n}$.

Consider $\sigma=x_{i_{1}}\ldots x_{i_{n}}$ a permutation of $S_{n}.\;$The
associated triangulation $\varphi(\sigma)$ of $P$ is defined by adding $n-1$
noncrossing diagonals to $P$ with the following algorithm:

\begin{enumerate}
\item  In $P^{(1)}=P$, join the two neighbors of $x_{i_{1}}$, i.e.
$\mathrm{pred}(x_{i_{1}})$ and $\mathrm{succ}(x_{i_{1}})$, by a diagonal.

%
%
%
%

\item  For each integer $2\leq k\leq n-1,$ consider the polygon $P^{(k)}$
obtained from $P^{(k-1)}$ by deleting the vertex $x_{i_{k-1}}$ and the edges
connected to it, then joining the two neighbors of $x_{i_{k}}$ in the polygon
$P^{(k)}$ by a diagonal.%

\begin{center}
\includegraphics[
height=7.5674cm,
width=10.1067cm
]%
{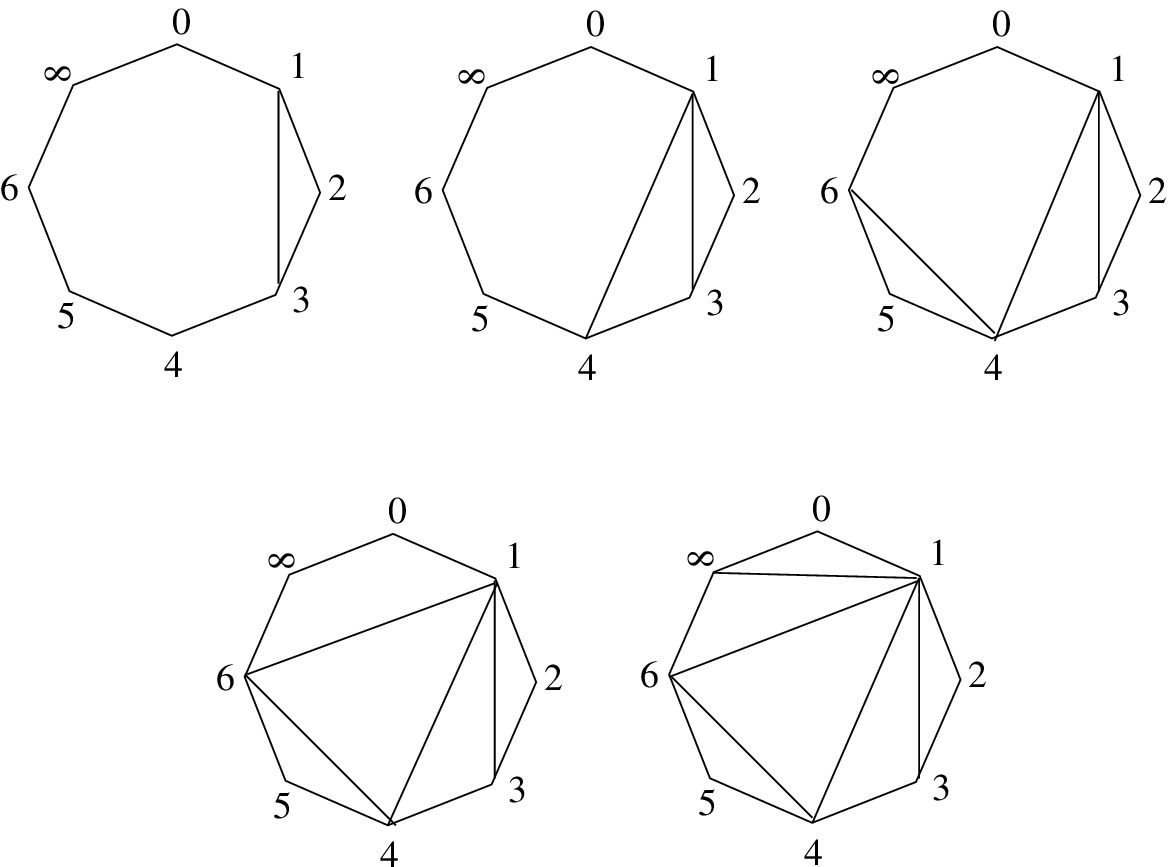}%
\\
Figure 4: The triangulation $\varphi(235461)$%
\label{fig4}%
\end{center}
\end{enumerate}

One easily verifies that the plane graph obtained when the procedure
terminates is a triangulation $T$ of $P.\;$

\bigskip

\noindent\textbf{Remark. }In \cite{Re}, the author uses a similar
combinatorial map, denoted $\eta$, from the symmetric group to the set of
triangulations.\ The definitions of the maps $\varphi$ and $\eta$ are quite
different, but one can verify that for any $\sigma=x_{i_{1}}\cdot\cdot\cdot
x_{i_{n}}\in S_{n},$ one has
\[
\eta(\sigma)=\varphi(x_{i_{n}}\cdot\cdot\cdot x_{i_{1}}).
\]

\bigskip

We shall now prove that each triangulation $T\in\mathcal{T}_{n}$ can be
represented by a suitable word $\sigma\in S_{n}$.

\begin{lemma}
\label{cutt}The map $\varphi:S_{n}\rightarrow\mathcal{T}_{n}$ is surjective.
\end{lemma}

\begin{proof}
We proceed by induction on $n.\;$ The lemma is trivial for $n=1$. Assume
$n\geq2$ and the statement true for $n-1$. Let $X=\{x_{1}<x_{2}<\ldots<x_{n}\}
\subset\mathbb{N}\setminus\{0\}$ and let $T$ be a triangulation of the
$(n+2)$-gon $P$ with vertex set $\widehat{X}=\{0,x_{1},x_{2},\ldots
,x_{n},\infty\}.$ Since $T$ contains at least two non-adjacent ears, some
vertex $x_{i} \in X$ must be an ear. Denote by $T^{\prime}$ the triangulation
obtained by \textit{cutting the ear $x_{i}$\;} in $T$, i.e. by deleting in $T$
the unique face containing the ear $x_{i}$. Then $T^{\prime}$ is a
triangulation of a convex $(n+1)$-gon $P^{\prime}$ on the vertex set
$\widehat{X}\setminus\{x_{i}\}$, so that $T^{\prime}\in\mathcal{T}_{n-1}$. By
the induction hypothesis, $T^{\prime}=\varphi(\sigma^{\prime})$ where
$\sigma^{\prime}=x_{i_{2}}\dots x_{i_{n}}$ is a permutation of the set
$X\setminus\{x_{i}\}$. Denoting $\sigma=x_{i} x_{i_{2}}\dots x_{i_{n}} \in
S_{n}$, we have $\varphi(\sigma)=T$ by construction.
\end{proof}

\bigskip

We shall refer to the procedure used in the above proof as the \textit{cutting
ear procedure}. The words in $S_{n}$ obtained with this procedure will be
called the \textit{readings} of the triangulation $T.\;$From the previous
lemma we have $\varphi(\sigma)=T$ for any reading $\sigma$ of $T.\;$The
\textit{canonical reading} of $T$ is the reading obtained by considering the
cutting ear procedure for which the vertex with the greatest label is deleted
at each step. This is equivalent to say that the readings of $T$ are the words
in $\varphi^{-1}(T)$ and the canonical reading of $T$ is the greatest reading
for the lexicographic order (see Example \ref{read}).

\subsection{The simple colored triangulation associated to a
word\label{suubsec_deffi}}

We shall now slightly generalize the preceding construction, by associating a
suitably colored triangulation to a word where letters are allowed to be
repeated. Consider a word $w$ of length $n$ on the alphabet $\mathcal{C}%
$.\ The \textit{evaluation} of $w$ is the $p$-uple $\mu=(\mu_{1},\ldots
,\mu_{p})\in\mathbb{N}^{p}$ where for any $k\in\{1,\ldots,p\},$ $\mu_{k}$ is
the number of occurrences of the color $c_{k}$ in $w.\;$We write for short
$\mathrm{eval(}w)=\mu.\;$ We denote by $\mathcal{C}_{n,\mu}$ the set of words
of length $n$ and evaluation $\mu$ on the alphabet $\mathcal{C}$. The
\textit{standardization} of $w$ on $X\;$ will be denoted by $\mathrm{std}(w)$.
Recall that $\mathrm{std}(w)$ is obtained by labelling from $x_{1}$ to
$x_{\mu_{1}}$ the occurrences of the color $c_{1}$ reading from left to right,
then from $x_{\mu_{1}+1}$ to $x_{\mu_{1}+\mu_{2}}$ the occurrences of $c_{2},$
and so on. (See the example below.)

For any $\sigma\in S_{n}$ considered as a word on $X=\{x_{1}<\cdot\cdot
\cdot<x_{n}\}$ of length $n,$ one associates the sequence $\Delta
(\sigma)=(\delta_{1},\ldots,\delta_{p})$ where $\delta_{1}$ is the longest
increasing sequence in $\sigma$ starting at $x_{1}$ of successive letters read
from left to right, $\delta_{2}$ the longest increasing sequence in $\sigma$
starting at $y_{1}=\max(\delta_{1})$ of successive letters and so on. Set
$L(\sigma)=\{\ell_{1},\ldots,\ell_{p}\}$, where $\ell_{k}$ is the length of
$\delta_{k}$ for any $k=1,\ldots,p,$ and set $S_{n}^{(\mu)}=\{\sigma\in
S_{n}\mid L(\sigma)=\mu\}$. Then, for any word $w \in\mathcal{C}_{n,\mu}$, its
standardization std$(w)$ belongs to $S_{n}^{(\mu)}$, and the standardization
map
\[
\mathrm{std}: \mathcal{C}_{n,\mu} \rightarrow S_{n}^{(\mu)}
\]
is a bijection.\ The corresponding inverse map is called the
\textit{destandardization} and denoted by $\mathrm{dstd}_{\mu}$.

\begin{example}
Suppose $X=\{1,\ldots,8\}$ and $\mathcal{C}=\{a,b,c,d\}.\;$For $w=bacbbacd$ we
have%
\[
\left\{
\begin{tabular}
[c]{c}%
$\mathrm{eval(}w)=(2,3,2,1),\text{ } \mathrm{std}(bacbbacd)=31645278,$\\
$\Delta(31645278)=(12,345,678),$ $L(\sigma)=(2,3,3)$.
\end{tabular}
\right.
\]
In particular $\sigma=31672485\in S_{n}^{(\mu)}$ and $\mathrm{dstd}_{\mu
}(\sigma)=baccabdb$.
\end{example}

To each $n$-tuple $\mu=(\mu_{1},\ldots,\mu_{p})\in\mathbb{N}^{p}$ we associate
the coloring%
\begin{equation}
\varepsilon_{\mu}=(\underset{\mu_{1}\text{ times}}{\underbrace{c_{1}%
,\ldots,c_{1}}},\underset{\mu_{2}\text{ times}}{\underbrace{c_{2},\ldots
,c_{2}}},\ldots,\underset{\mu_{p}\text{ times}}{\underbrace{c_{p},\ldots
,c_{p}}})\in\mathbb{N}^{n}. \label{epsilon(mu)}%
\end{equation}
Our purpose now is to define a map
\[
\Phi:\mathcal{C}_{n,\mu}\rightarrow\mathcal{ST}_{n}(\varepsilon_{\mu}),
\]
where $\mathcal{ST}_{n}(\varepsilon_{\mu})$ is the set of simple colored
triangulations of $P$ with coloring $\varepsilon_{\mu}$. To the word
$w\in\mathcal{C}_{n,\mu}$ we associate the colored triangulation
$\Phi(w)=T_{\varepsilon_{\mu}}$, where $T=\varphi(\mathrm{std}(w))$.

\begin{lemma}
Let $\sigma$ be a word in $S_{n}$ on the alphabet $X = \{x_{1} <\ldots<
x_{n}\}$ and let $T=\varphi(\sigma)$ be the associated triangulation of the
polygon $P$ with vertex set $\widehat{X}=\{0 < x_{1} <\ldots< x_{n} <
\infty\}$. For each $i=1,\ldots, n-1$, denote by $t_{i} \in\widehat{X}$ the
third vertex of the unique face of~ $T$ containing the edge $\{x_{i},
x_{i+1}\}$. Then \ 

\begin{enumerate}
\item  The letter $x_{i}$ is on the left of the letter $x_{i+1}$ in $\sigma$
if and only if $t_{i}<x_{i}.$

\item  If $\sigma$ contains the increasing sequence $x_{i},x_{i+1}%
,\ldots,x_{j}$ with $j \geq i+2$ in its left to right reading, then $t_{j-1}
\leq t_{j-2} \leq\ldots\leq t_{i} < x_{i}.$

\item  For any word $w$ on $\mathcal{C}$, $\Phi(w)$ is a simple colored triangulation.
\end{enumerate}
\end{lemma}

\begin{proof}
\ \ 

\noindent$(1)$ Suppose that $x_{i}$ is on the left of the letter $x_{i+1}$ in
$\sigma=y_{1}y_{2}\cdots y_{n}$. Assume $x_{i}=y_{k}$. Then, when constructing
$T=\varphi(\sigma)$ step-by-step, the vertex $x_{i}$ is an ear in the
subtriangulation $T^{\prime}=\varphi({y_{k}\cdots y_{n}})$, and its successor
in $Y=\widehat{X}\setminus{\{y_{1},\ldots,y_{k-1}\}}$ is $x_{i+1}.$ Since the
unique face of $T^{\prime}$ containing $x_{i}$ is $\{\mathrm{pred}_{Y}%
(x_{i}),x_{i},x_{i+1}\}$, we must have $t_{i}=\mathrm{pred}_{Y}(x_{i}%
)<x_{i}.\;$Conversely, assume that $x_{i}$ is on the right of $x_{i+1}$ in
$\sigma$. Then $x_{i+1}$ is an ear in some subtriangulation $T^{\prime}$ of
$T$, and its predecessor in the vertex set $Y$ of $T^{\prime}$ must be $x_{i}%
$. Thus $t_{i}$ is the successor of $x_{i+1}$ in $Y$, and therefore
$x_{i}<x_{i+1}<t_{i}$. This proves assertion $(1)$.

\smallskip

\noindent$(2)$ We know by $(1)$ that $t_{i} < x_{i}$. By reasoning inductively
on $j$, it suffices to show that $t_{i+1} \leq t_{i}$. As above, the vertex
$x_{i}$ is an ear in some subtriangulation $T^{\prime}$ of $T$ with vertex set
$Y \subset\widehat{X}$, and $t_{i}$ is the predecessor of $x_{i}$ in $Y$.
Since $t_{i+1}$ is the predecessor of $x_{i+1}$ in some subset of $Y$, and
since $t_{i+1}\not =x_{i}$ (as there is an edge $\{x_{i+2}, t_{i+1}\}$ but no
edge $\{x_{i+2}, x_{i}\}$), we must have $t_{i+1} \leq t_{i}$ as claimed.

\smallskip

\noindent$(3)$ Note first that by definition of $\Phi,$ $\Phi(w)$ is an
increasing colored triangulation.\ Set $\sigma=\mathrm{std}(w)$ and
$T=\varphi(\sigma)$.$\;$Consider a face $F$ of $\Phi(w)$ with two vertices
colored by $c.\;$In $T$, $F$ has two vertices labelled by $x_{r}$ and $x_{s}$
such that $x_{r}<x_{s}$.\ All the letters $y\in X$ with $x_{r}\leq y\leq
x_{s}$ label vertices colored by $c$ in $\Phi(w).\;$Thus they correspond to
letters $c$ in $w$. Since $\sigma=\mathrm{std}(w)$, they must appear in
increasing order in the left to right reading of $\sigma.$ Suppose $r-s>1$,
that is $x_{r}x_{s}$ is an inner diagonal of $T$.\ Then $x_{r+1}$ appears in
$\sigma$ on the right of $x_{r}.\;$When $\varphi(\sigma)$ is constructed, a
diagonal joining the vertex labelled $x_{r+1}$ to a vertex labelled by $t<$
$x_{r+1}$ is drawn.\ Since there are no intersections between the inner
diagonals of $\varphi(\sigma),$ this gives a contradiction for $F$ is a face
with an edge joining $x_{r}$ and $x_{s}$ and $s>r+1.\;$This means that
$x_{s}=x_{r+1}$.\ By $(1),$ we obtain that the third face of $F$ is labelled
by $t_{r}<x_{r}.$ Thus $\Phi(w)$ is a simple colored triangulation.
\end{proof}

\bigskip%

\begin{center}
\includegraphics[
height=7.0094cm,
width=12.1166cm
]%
{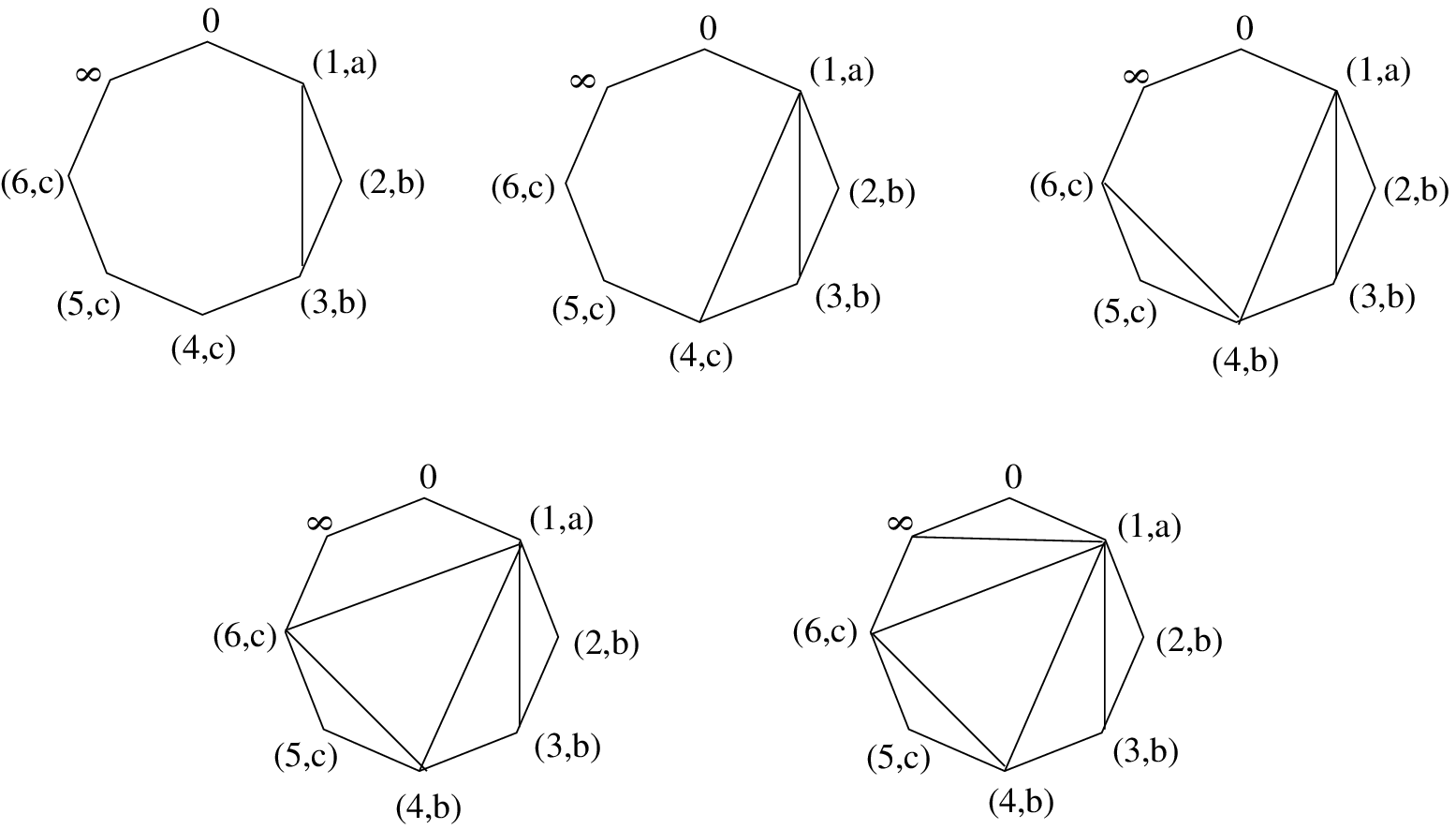}%
\\
Figure 5: The simple triangulation $\Phi(bbcbca)$%
\label{fig5}%
\end{center}

Consider $T_{\varepsilon_{\mu}}\in\mathcal{ST}_{n}(\varepsilon_{\mu})$ (see
(\ref{epsilon(mu)})).$\;$The readings of the simple colored triangulation
$T_{\varepsilon_{\mu}}$ are the words\ obtained by applying cutting ear
procedures on $T$ and by forming the words of $\mathcal{C}_{n,\mu}$ obtained
from the successive colors of the ears deleted instead of their labels. Since
the cutting ear procedure depends only of the inner diagonals of the
triangulation $T$ and not on the labels or colorings of its faces, there is a
one-to-one correspondence between the readings of $T_{\varepsilon_{\mu}}$ and
those of $T$. Moreover by definition of $\Phi,$ the readings of $T$ coincide
with the standardized of the readings of $T_{\varepsilon_{\mu}}.$ The
canonical reading of $T_{\varepsilon_{\mu}}$ is the greatest reading of
$T_{\varepsilon_{\mu}}$ for the lexicographic order. By proceeding as in Lemma
\ref{cutt}, we prove that for any reading $w$ of $T_{\varepsilon_{\mu}},$ we
have $\Phi(w)=T_{\varepsilon_{\mu}}.$ In particular the map $\Phi$ is surjective.

\begin{example}
\label{read}The readings of the triangulation obtained in Figure 4 are the
permutations $\{235461,253461,523461\}$.\ Thus the readings of the colored
triangulation obtained in the above figure are the words
$\{bbcbca,bcbbca,cbbbca\}$.
\end{example}

\bigskip

The map $\Phi$ can alternatively be thought of as an insertion scheme for the
simple colored triangulations.\ Indeed consider a simple colored triangulation
$T_{\varepsilon}$.\ Suppose that the vertices of $T$ are indexed by the
integers of the set $\widehat{X}=\{0<x_{1}<\cdot\cdot\cdot<x_{n}<\infty\}$.
Let $x$ be a positive integer such that $x\notin X$ and $x_{1}<x<x_{n}$.
Consider the pair $(x,c)$ where $c\in\mathcal{C}.$ There exists a unique
vertex $v$ in $T_{\varepsilon}$ which is not colored by $c$ (thus is colored
by $b<c$ since $\varepsilon$ is an increasing coloring) and such that
$v^{\prime}=\mathrm{succ}(v)$ is colored by $c.$ The insertion of $(x,c)$ in
$T_{\varepsilon}$ (denoted $(x,c)\rightarrow T_{\varepsilon}$ for short) is
defined as follows:

\begin{itemize}
\item  Add a vertex $v^{\ast}$ indexed by $(x,c)$ between the vertices $v$ and
$v^{\prime}$ in $T_{\varepsilon}.$

\item  Draw the two edges joining $v^{\ast}$ to $v$ and $v^{\ast}$ to
$v^{\prime}$.

\item  Color the face defined by $v,v^{\ast}$ and $v^{\prime}$ with $c.$
\end{itemize}

\begin{proposition}
\ \ \ 

\begin{enumerate}
\item  The output of the insertion $(x,c)\rightarrow T_{\varepsilon}$ is a
simple colored triangulation.

\item  For any word $w=c_{1}\cdot\cdot\cdot c_{n}$ on $\mathcal{C}_{n,\mu}$,
we have
\begin{equation}
\Phi(w)=(x_{1},c_{1})\rightarrow(x_{2},c_{2})\rightarrow\cdot\cdot
\cdot\rightarrow(x_{n},c_{n})\rightarrow T_{0} \label{ins}%
\end{equation}
where $\sigma=\mathrm{sdt}(w)$ and $T_{0}$ is the triangulation with one edge
joining the vertices $0$ and $\infty.$
\end{enumerate}
\end{proposition}

\begin{proof}
\ 

\noindent$(1)$ If $F$ is a face of $(x,c)\rightarrow T_{\varepsilon}$ with two
vertices colored by the same color $a,$ $F$ is face of $T_{\varepsilon}$ or
$F=vv^{\ast}v^{\prime}.\;$In the first case, the third vertex of $F$ is
labelled by $0$ or colored by $a^{\prime}<a$ since $T_{\varepsilon}^{\prime}$
is simple. In the second case, the third vertex of $F$ is $v^{\prime}$ which
is labelled by $0$ or colored by $b^{\prime}<c.$

\noindent$(2)$ We proceed by induction on $n.\;$When $n=1,$ the assertion is
true since $\Phi(w)$ has only one face. Suppose that (\ref{ins}) holds for
$n-1$ and consider $w\in\mathcal{C}_{n,\mu}.\;$Write $w=t_{1}w$ where
$w^{\prime}=c_{2}\cdot\cdot\cdot c_{n}$. Then we have by induction
$\Phi(w^{\prime})=(x_{2},c_{2})\rightarrow\cdot\cdot\cdot\rightarrow
(x_{n},c_{n})\rightarrow T_{0}.$ Thus it suffices to show that $(x_{1}%
,c_{1})\rightarrow\Phi(\sigma^{\prime})=\Phi(\sigma)$ which immediately
follows from the definition of $\Phi$ and the description of the insertion
algorithm given above.
\end{proof}

\bigskip

\noindent\textbf{Remark. } Binary trees and triangulations are known to be
formally equivalent combinatorial objects.\ This analogy can be extended to
binary search trees which corresponds to simple colored triangulations. The
previous insertion scheme can be then regarded as an analogue of Knuth's
insertion algorithm on binary search trees \cite{KN}.%

\begin{center}
\includegraphics[
trim=0.000000in 0.000000in 0.000000in 0.005526in,
height=13.5817cm,
width=10.1067cm
]%
{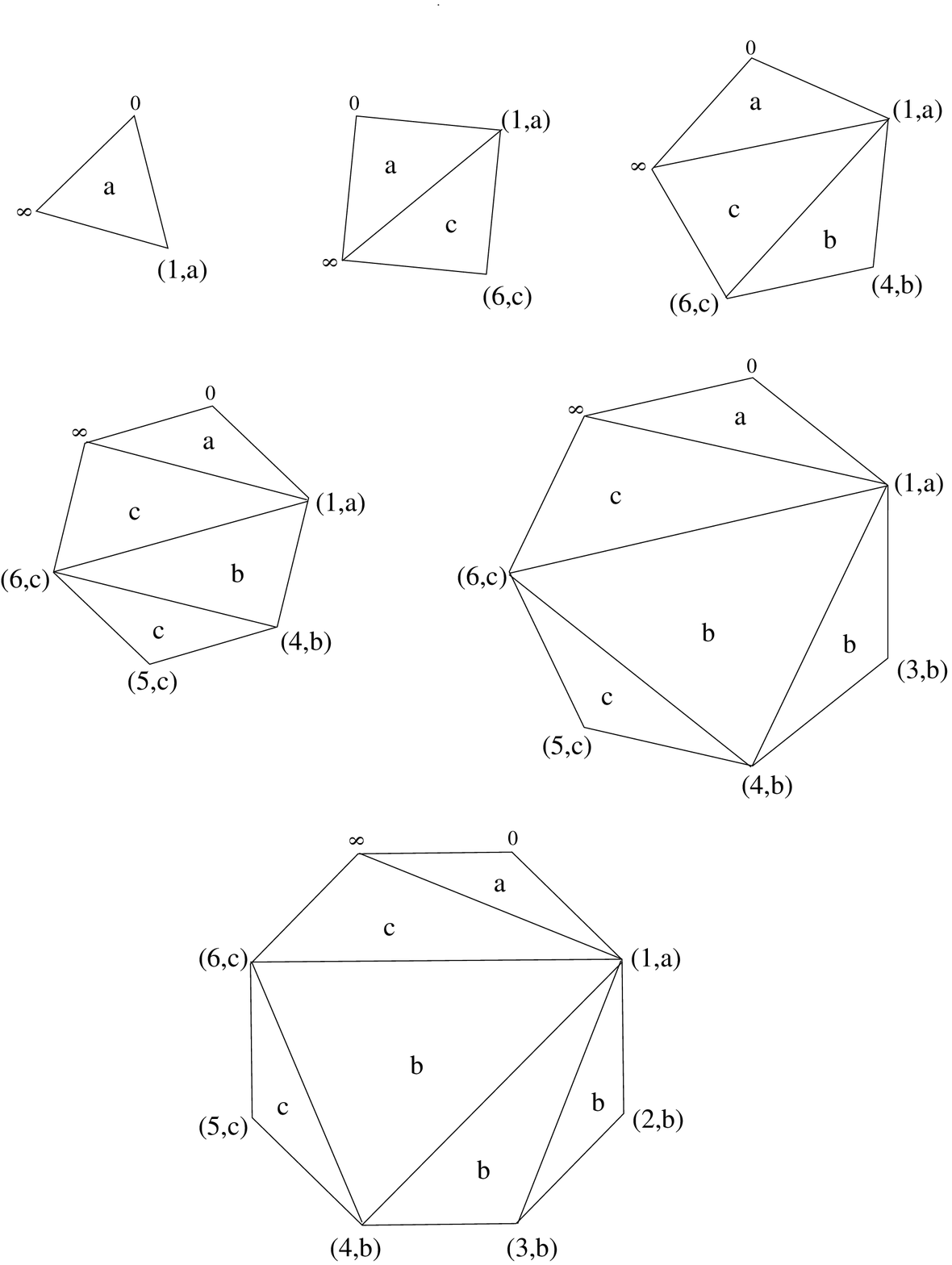}%
\\
Figure 6: $\Phi(w)$ with $w=bbcbca$ and $\sigma=235461$%
\label{fig10}%
\end{center}

\subsection{The Cayley graph of $S_{n}$ and the flip graph $\mathcal{F}_{n}$}

In this paragraph we prove that the above map $\varphi:S_{n}\rightarrow
\mathcal{T}_{n}$ is in fact a surjective \textit{morphism} from the Cayley
graph of $S_{n},$ $\mathrm{Cay}(S_{n}),$ to the flip graph $\mathcal{F}_{n}.$
Here $\mathrm{Cay}(S_{n})$ is the graph on the vertex set $S_{n}$, with edges
all pairs $\{\sigma,\sigma^{\prime}\}$ such that $\sigma,\sigma^{\prime}\in
S_{n}$ differ by an elementary transposition. Since we have identified the
permutations of $X=\{x_{1}<\ldots<x_{n}\}$ with the standard words of length
$n$ with letters in $X,$ this is equivalent to say that $\sigma^{\prime}$ is
obtained by switching two consecutive letters of $\sigma$, that is
\begin{equation}%
\begin{array}
[c]{lcl}%
\sigma & = & uxzv\\
\sigma^{\prime} & = & uzxv\label{edge}%
\end{array}
\end{equation}
where $x,z$ are letters in $X$ and $u,v$ standard words on $X.$ By symmetry,
we can suppose $x<z$.

\begin{lemma}
\label{lemfond} Let $\sigma=uxzv\text{ and }\sigma^{\prime}=uzxv$, where as
above $x < z \in X$ and $u,v$ are words on $X$.

\begin{enumerate}
\item  Suppose that $v$ does not contain any letter $y\in X$ such that
$x<y<z$. Then the associated triangulations $\varphi(\sigma)$ and
$\varphi(\sigma^{\prime})$ differ by a diagonal flip. In other words, the edge
$\{\sigma, \sigma^{\prime}\}$ in $\mathrm{Cay}(S_{n})$ is mapped by $\varphi$
to an edge in the flip graph $\mathcal{F}_{n}.$

\item  Suppose that $v$ contains a letter $y$ such that $x<y<z$.\ Then
$\varphi(\sigma)=\varphi(\sigma^{\prime}).$ Here, the edge $\{\sigma,
\sigma^{\prime}\}$ in $\mathrm{Cay}(S_{n})$ is contracted by $\varphi$ to a
single vertex in $\mathcal{F}_{n}.$
\end{enumerate}
\end{lemma}

\begin{proof}
\ 

\noindent$(1)$ Since $\sigma$ is a standard word, the letters $y$ such that
$x<y<z$, if any, belong to $u.\;$Set $Y=X-\{u\}.\;$Then $\mathrm{pred}%
_{\widehat{Y}}(z)=x.\;$Set $x^{\prime}=\mathrm{pred}_{\widehat{Y}}(x)$ and
$z^{\prime}=\mathrm{succ}_{\widehat{Y}}(z).$ Then $\varphi(uzx)$ is obtained
from $\varphi(uxz)$ by flipping the diagonal $x^{\prime}z$ in the
quadrilateral $x^{\prime}xzz^{\prime}.\;$Thus $\varphi(\sigma)$ and
$\varphi(\sigma^{\prime})$ also differ by a diagonal flip in $\mathcal{F}_{n}.$

\noindent$(2)$ This times there exists a letter $y$ such that $x<y<z$ in
$Y.\;$Hence $\mathrm{pred}_{\widehat{Y}}(z)\neq x$.\ This implies immediately
that $\varphi(uzx)=\varphi(uxz)$ and thus $\varphi(\sigma)=\varphi
(\sigma^{\prime}).$
\end{proof}

\begin{proposition}
The map $\varphi$ is a surjective morphism of graphs from $\mathrm{Cay}%
(S_{n})$ to $\mathcal{F}_{n}.$
\end{proposition}

\begin{proof}
We have already obtained that $\varphi$ is a surjective map.\ Now suppose that
the permutation $\sigma$ and $\sigma^{\prime}$ belongs to the same edge in
$\mathrm{Cay}(S_{n}).\;$Then they can be written as in (\ref{edge}).\ The
previous lemma implies that we will either have $\varphi(\sigma)=\varphi
(\sigma^{\prime})$ or the pair $\varphi(\sigma), \varphi(\sigma^{\prime})$
differs by exactly one diagonal flip.\ Thus $\varphi$ is a morphism of graphs,
as claimed.
\end{proof}

\subsection{Sylvester relations}

In this section we shall determine the fibers of the morphism $\varphi
:S_{n}\rightarrow\mathcal{T}_{n}$. It turns out that they can be characterized
in terms of the \textit{sylvester\footnote{By analogy with the french word
\textit{sylvestre} which means \textit{forestal} in english} relations}
introduced by Hivert, Novelli and Thibon in the context of binary search trees
\cite{HNT}.\ Two words $w_{1}$ and $w_{2}$ on a totally ordered alphabet
$\mathcal{A}$ are said to be \textit{sylvester adjacent} if there exist three
words $u,u^{\prime},u^{\prime\prime} \in\mathcal{A}^{\ast}$ and three letters
$x\leq y<z \in\mathcal{A}$ such that
%
%
%
%
\begin{equation}%
\begin{array}
[c]{rcl}%
w_{1} & = & uxzu^{\prime}yu^{\prime\prime}\\
w_{2} & = & uzxu^{\prime}yu^{\prime\prime}. \label{sylrel}%
\end{array}
\end{equation}
Two words $w$ and $w^{\prime}$ are \textit{sylvester congruent} if there
exists a chain of words
\[
w=w_{1},w_{2},\ldots,w_{k}=w^{\prime}%
\]
such that for any $i=1,\ldots,k-1,$ the words $w_{i}$ and $w_{i+1}$ are
sylvester adjacent. The sylvester congruence is a congruence denoted by
$\equiv$ on the free monoid $\mathcal{A}^{\ast}$ on $\mathcal{A}$.\ We write
$[w]$ for the sylvester class of the word $w \in\mathcal{A}^{\ast}.$

The sylvester congruence has the following remarkable property.

\begin{lemma}
\ \label{lem-stand}

\begin{enumerate}
\item  Let $w_{1}$ and $w_{2}$ two words on $\mathcal{A}$ such that
$w_{1}\equiv w_{2}.\;$Then we have $\mathrm{std}(w_{1})\equiv\mathrm{std}%
(w_{2}).$

\item  Conversely if $w_{1}$ and $w_{2}$ have the same evaluation $\mu$ and
verify $\mathrm{std}(w_{1})\equiv\mathrm{std}(w_{2}),$ we have $w_{1}\equiv
w_{2}.$
\end{enumerate}
\end{lemma}

\begin{proof}
The statement follows immediately from the definition of the standardization
map and the definition (\ref{sylrel}) of sylvester adjacency.
\end{proof}

\bigskip

We now turn to the description of the fibers of $\varphi:S_{n}\rightarrow
\mathcal{T}_{n}$ in terms of the sylvester relations.

\begin{proposition}
\label{prop_sylvSn}Consider $\sigma_{1}$ and $\sigma_{2}$ two permutations in
$S_{n}.\;$Then
\[
\varphi(\sigma_{1})=\varphi(\sigma_{2})\Longleftrightarrow\sigma_{1}%
\equiv\sigma_{2}.
\]
In other words, two permutations are sylvester congruent if and only if they
are readings of the same triangulation.
\end{proposition}

\begin{proof}
The proposition is immediate for $n\leq2$, and so we assume now $n\geq3$.

Suppose that $\sigma_{1}\equiv\sigma_{2}.\;$By induction, it is enough to
prove that $\varphi(\sigma_{1})=\varphi(\sigma_{2})$ when $\sigma_{1}$ and
$\sigma_{2}$ are sylvester adjacent. This directly follows from $(2)$ of Lemma
\ref{lemfond}. This shows in particular that for any triangulation $T,$ the
fiber $\varphi^{-1}(T)$ is a nonempty disjoint union of sylvester classes.

To obtain the left to right part of the proposition, observe first that all
the words belonging to the same sylvester class have the same last letter.\ By
an easy induction, this implies that each sylvester class contains words $w$
which verify the separation property
\begin{equation}
w=uvy \label{sp}%
\end{equation}
where $y$ is a letter and $u,v$ are words either empty or containing letters
greater than $y$ and smaller than $y$, respectively. Such a word $w$ is not
unique. To obtain a normal form it suffices to impose that $u$ and $v$ are
empty or verify themselves the separation property. Denote by $f_{k}$ the
number of such normal forms for the permutations of $S_{k}.\;$Then the
separation property (\ref{sp}) gives the following recurrence formula for the
numbers $f_{k}$:%
\[
f_{n}=\sum_{i=1}^{n}f_{i-1}f_{n-i}.
\]
Moreover we have $f_{0}=1$ and $f_{1}=1$.\ Thus the above recurrence formula
is the recurrence formula for the Catalan numbers $c_{n}$.\ This shows that
the number of sylvester classes is equal to $c_{n}.$ Since $\varphi$ is
surjective, there are exactly $c_{n}$ nonempty fibers $\varphi^{-1}(T)$ which
all contains at least a sylvester class. This means that each fiber is a
sylvester class and immediately yields the desired implication $\varphi
(\sigma_{1})=\varphi(\sigma_{2})\Longrightarrow\sigma_{1}\equiv\sigma_{2}.$
\end{proof}

\bigskip

The previous result easily extends to the setting of colored triangulations.

\begin{theorem}
\label{Th_sylvT}Consider $\mu=(\mu_{1},\ldots,\mu_{p})\in\mathbb{N}^{p}$ such
that $\mu_{1}+\cdot\cdot\cdot+\mu_{p}=n$ and $w_{1},w_{2}$ two words of
evaluation $\mu$ on the alphabet $\mathcal{C}$.\ Then
\[
\Phi(w_{1})=\Phi(w_{2})\Longleftrightarrow w_{1}\equiv w_{2}.
\]
Two words on $\mathcal{C}$ are sylvester congruent if and only if they are
readings of the same simple colored triangulation.
\end{theorem}

\begin{proof}
Since $w_{1}$ and $w_{2}$ have the same evaluation $\mu,$ we deduce from Lemma
\ref{lem-stand} and Proposition \ref{prop_sylvSn} the equivalencies:
\[
\Phi(w_{1})=\Phi(w_{2})\Longleftrightarrow\varphi(\mathrm{std}(w_{1}%
))=\varphi(\mathrm{std}(w_{2}))\Longleftrightarrow\mathrm{std}(w_{1}%
)\equiv\mathrm{std}(w_{2})\Longleftrightarrow w_{1}\equiv w_{2}.
\]
\end{proof}

\bigskip

\noindent\textbf{Remarks.}

\noindent$\mathrm{(i)}$ The interpretation of the sylvester congruence we use
in this paper is not that originally given in \cite{HNT} where the sylvester
classes are defined as the fibers of the map which associated to each word the
binary search tree obtained via Knuth's insertion algorithm on binary trees.
This means that we have chosen to parametrize the sylvester classes by simple
colored triangulations rather than binary search trees.\ To obtain the binary
tree corresponding to a simple signed triangulation of canonical reading $w$,
one applies the insertion algorithm on binary trees starting from
$w.\;$Conversely, the triangulation associated to a binary tree with right to
left postfix reading $w$, is the triangulation $\Phi(w)$. Note that $w$ is
then the canonical reading of $\Phi(w)$.

\noindent$\mathrm{(ii)}$ It is very easy to obtain the sylvester class of a
word $w$ from its associated colored triangulation $\Phi(w).\;$Indeed, $[w]$
is simply the set of readings of $\Phi(w)$ defined in \ref{suubsec_deffi}.

\noindent$\mathrm{(iii)}$ For any simple triangulation $T_{\varepsilon},$ we
know by Theorem \ref{Th_sylvT} that the set $[T_{\varepsilon}]$ of all the
different readings of $T_{\varepsilon}$ is a plactic class. The same property
holds for the set of readings $[T]$ of $T.\;$Then the map $\mathrm{std}$ is a
one-to-one correspondence between $[T_{\varepsilon}]$ and $[T]$.

\begin{proposition}
\label{prop_flip}Two triangulations $T_{1}$ and $T_{2}$ differ by a diagonal
flip if and only if there exist a reading $w_{1}$ of $T_{1}$ and a reading
$w_{2}$ of $T_{2}$ such that
%
%
%
%
\begin{equation}%
\begin{array}
[c]{rcl}%
w_{1} & = & uxzv\\
w_{2} & = & uzxv \label{cond}%
\end{array}
\end{equation}
where $x, z$ belong to $X$ and where $u, v$ are words such that $v$ contains
no letter $y\in X$ satisfying $x<y<z$ or $z<y<x.$
\end{proposition}

\begin{proof}
If $w_{1}$ and $w_{2}$ are readings respectively of $T_{1}$ and $T_{2}$
verifying (\ref{cond}), they differ by a diagonal flip by $(1)$ of Lemma
\ref{lemfond}.\ Conversely suppose that $T_{1}$ and $T_{2}$ differ by a
diagonal flip.\ Denote by $Q=x^{\prime}xzz^{\prime}$ the quadrilateral in
which the diagonal flip happens.\ One can suppose that $x^{\prime}xzz^{\prime
}$ is the clockwise reading of the vertices of $Q$ and consider that the
diagonal $x^{\prime}z$ is flipped in $T_{1}$ to give the diagonal $xz^{\prime
}$ in $T_{2}.$ Consider the triangulations $K,L$ and $M$ whose vertices are
respectively indexed by the integers $k,l$ and $m$ in $\widehat{X}$ such that
\[
x^{\prime}\leq k\leq x,\text{ }x\leq l\leq z\text{ and }z\leq m\leq z^{\prime
}.
\]
In $T_{1}$ and $T_{2},$ the vertices $x^{\prime}$ and $x$ are connected to
edges which are not edges of $K.\;$We have the same property for the vertices
$z$ and $x$ in $L$, and for the vectors $z$ and $z^{\prime}$ in $M.\;$Thus
there is at least a reading of $T_{1}$ and $T_{2}$ which starts with the
vertices of $K,L,M$ distinct of $x^{\prime},x,z$ and $z^{\prime}.\;$Choose one
of these readings and denote by $u$ the word obtained.\ Then $u$ contains all
the vertices labelled by an integer $y$ such that $x<y<z.$ Moreover, we can go
on the lecture of $T_{1}$ by reading successively $x,z$ and next choose a
reading $v$ of the remaining triangulation.\ Similarly, we can go on the
lecture of $T_{2}$ by reading successively the vertices $z,x$ and next form
the word $v$ as previously. Finally the readings we obtain for $T_{1}$ and
$T_{2}$ are respectively $uxzv$ and $uzxv$ and the integers $y$ such that
$x<y<z$ are in $u.$
\end{proof}

\bigskip

\noindent\textbf{Remark. } This proof shows also that the letters $x$ and $z$
appearing in (\ref{cond}) are the labels of the faces of the quadrilateral of
$T_{1}$ in which the flip happens.

\section{Reformulations of the four color theorem}

It is well known that, considered as a statement of graph theory, the four
color theorem is equivalent to say that every simple finite planar graph
admits a proper four coloring of its vertices. In \cite{El1}, the first author
obtains a reformulation of this theorem in terms of signed paths between
triangulations of polygons.\ Let us recall briefly the main ideas which permit
this reformulation. Tutte has proved in \cite{tu} that every $4$-connected
finite planar graph is hamiltonian (i.e. admits a cycle which visits each
vertex exactly once). This implies that it suffices to prove the four color
theorem for hamiltonian planar triangulations of the sphere $S^{2}$.\ For such
a triangulation $\frak{T}$ , it is natural to confine the hamiltonian path $P$
on the equator of $S^{2}$.\ The path $P$ can be then regarded as a polygon and
one defines two triangulations $T_{1}$ and $T_{2}$ of $P$ by considering the
sub-graphs of $\frak{T}$ lying respectively in the northern and southern
hemispheres of $S^{2}$.\ Conversely, two triangulations $T_{1}$ and $T_{2}$ of
the same polygon $P$ define a triangulation $\frak{T}=T_{1}\cup T_{2}$ of
$S^{2}$ obtained by embedding them in the two hemispheres of $S^{2}$ and
gluing them along their common boundary $P$ confined on the equator. We denote
by $F(\frak{T})$ the set of faces of $\frak{T}$ and by $F_{v}$ the subset of
faces incident to some vertex $v$.\ A signing of $\frak{T}$ is a map
$\varepsilon:F(\frak{T)}\rightarrow\{\pm1\}$ which associates to each face of
$\frak{T}$ one of the integers $1$ or $-1.\;$We then denote by $\frak{T}%
_{\varepsilon}$ the signed triangulation obtained.\ Given $v$ a vertex of
$\frak{T}_{\varepsilon}$ write%
\[
s_{\varepsilon}(v)=\sum_{F\in F(\frak{T})}\varepsilon(F)
\]
for the sum of signs of the faces $F$ incident to $v.\;$The signing
$\varepsilon$ is a \textit{Heawood signing} if at each vertex $v$ of
$\frak{T}_{\varepsilon}\frak{,}$ one has%
\[
s_{\varepsilon}(v)\equiv0\operatorname{mod}3\text{.}%
\]
We will then say that $\frak{T}_{\varepsilon}$ is a Heawood signed
triangulation.\ At the end of the $19$-th century, Heawood \cite{hea} has
proved that $\frak{T}$ admits a proper four coloring of its vertices if and
only if there exists a Heawood signing on its faces. Now consider a
triangulation $T$ of the polygon $P.\;$As observed in \cite{tu} there is a
very simple way to obtain from $T$ a triangulation of the sphere which admits
a Heawood signing. Indeed, for $\frak{T}=T\cup T,$ every northern face has a
corresponding southern face having the same three vertices.\ So its suffices
to define $\varepsilon$ such that these faces have opposite signs to obtain a
Heawood signing. The diagonal flip operations on the triangulation $\frak{T}$
are natural geometrical transformations which yields new triangulations from
$\frak{T}$.\ Two adjacent faces $F$ and $F^{\prime}$ in the same plane
trinagulation $T_{1}$ or $T_{2}$ defining $\frak{T}$ form a quadrilateral $Q$
and the bound between $F$ and $F^{\prime}$ coincide with a diagonal $D$ of
$Q.\;$The diagonal flip operation in $Q$ delete the diagonal $D$ and replace
it by the opposite diagonal of $Q.$ To obtain diagonal flips operations on
Heawood signed triangulations (that is which preserve the Heawood property),
one has to restrict the authorized diagonals flips to what we call
\textit{signed flips}, defined as flips for which the signs of the two faces
of $Q$ are the same and are changed into their opposite during the flip
operation. Consider $\frak{T}_{\varepsilon}$ a hamiltonian planar
triangulation of the sphere and denote by $T_{1,\varepsilon_{1}}$ and
$T_{2,\varepsilon_{1}}$ the two signed plane triangulations such that
$\frak{T}_{\varepsilon}=T_{1,\varepsilon_{1}}\cup T_{2,\varepsilon_{1}}%
$.\ Suppose that there exists a sequence of signed flips from
$T_{1,\varepsilon_{1}}$ to $T_{2,\varepsilon_{1}}.\;$Since $\frak{T}%
_{(\varepsilon_{1},-\varepsilon_{1})}=T_{1,\varepsilon_{1}}\cup
T_{1,-\varepsilon_{1}}$ is a Heawood signed triangulation and by using that
the signed flips preserve the Heawood property, one then obtains that
$\frak{T}_{\varepsilon}$ is a Heawood triangulation of the sphere. Thus the
existence of a signed flip sequence between any two triangulations of a
polygon implies the four color theorem.\ This is the result obtained in
\cite{El1}.\ The converse is true as proved by Gravier and Payan in \cite{GP}.

\subsection{Signed flips and the four color theorem\label{subsec_fourcol}}

Consider a triangulation $T$ of a convex $(n+2)$-gon $P$. For any
$\varepsilon\in\{-,+\}^{n}$, denote by $T_{\varepsilon}$ the signed
triangulation obtaining by signing $T$ following $\varepsilon.\;$Note that
$T_{\varepsilon}$ is not simple in general.\ Let $\Sigma(T_{\varepsilon})$ be
the set of all signed triangulations obtained by applying a sequence of signed
flips starting from $T_{\varepsilon}$. For any signed triangulation
$U_{\varepsilon^{\prime}}$ belonging to $\Sigma(T_{\varepsilon})$, we will say
that there exists a \textit{signed path} between $T_{\varepsilon}$ and
$U_{\varepsilon^{\prime}}$.

\begin{lemma}
\label{lem-esse}Suppose that $\Sigma(T_{\varepsilon})$ contains signed
triangulations $U_{\varepsilon^{\prime}}$ and $U_{\varepsilon^{\prime\prime}}$
with the same underlying triangulation $U$.\ Then $\varepsilon^{\prime
}=\varepsilon^{\prime\prime}$.
\end{lemma}

\begin{proof}
Since the signed sphere triangulation $T_{\varepsilon}\cup_{P} T_{-\varepsilon
}$ has the Heawood property, and since signed flips preserve this property, it
follows that the signed sphere triangulation $U_{\varepsilon^{\prime}}\cup_{P}
U_{-\varepsilon^{\prime\prime}}$ also has the Heawood property. We now deduce
from this that $\varepsilon^{\prime}=\varepsilon^{\prime\prime}$ by induction
on the number $n+2$ of vertices. The statement is trivial for $n=1$, as there
is only one face in $U$. Assume $n \geq2$, and let $v$ be an ear in $U$. Thus
$v$ is contained in a unique face $F$ of $U$, and therefore is contained in
exactly two faces of $U \cup_{P} U$, namely one copy of $F$ on each
hemisphere. Since the signs of these two faces must sum up to 0 mod 3 in
$U_{\varepsilon^{\prime}}\cup_{P} U_{-\varepsilon^{\prime\prime}}$ by the
Heawood property, it follows that $\varepsilon^{\prime}(F)=\varepsilon
^{\prime\prime}(F)$. Let $Q$ denote the polygon obtained by contracting one of
the two edges of $P$ containing $v$, let $V$ denote the triangulation of $Q$
obtained by cutting the ear $v$ in $U$, and let $\mu^{\prime}$, respectively
$\mu^{\prime\prime}$, denote the restrictions of $\varepsilon^{\prime}$ and
$\varepsilon^{\prime\prime}$ to the faces of $V$. Then the signed sphere
triangulation $V_{\mu^{\prime}}\cup_{Q} V_{-\mu^{\prime\prime}}$ still has the
Heawood property, as easily seen. It follows by the induction hypothesis that
$\mu^{\prime}=\mu^{\prime\prime}$. Therefore, $\varepsilon^{\prime
}=\varepsilon^{\prime\prime}$ as claimed.
\end{proof}

\bigskip

\noindent\textbf{Remark. }By the previous lemma, there exists a signed path
between the signed triangulations $T_{\varepsilon}$ and $U_{\varepsilon
^{\prime}}$ only if there exists a path without loop between their underlying
triangulations $T$ and $U$ in $\mathcal{F}_{n}$.

\bigskip

Write $\left|  \Sigma(T_{\varepsilon})\right|  $ for the set of triangulations
obtained by deleting the signs $-$ and $+$ in the signed triangulations of
$\Sigma(T_{\varepsilon}).\;$We deduce from \cite{El1} and \cite{GP} the
following reformulation of the four color theorem.

\begin{theorem}
\label{ref1}The four color theorem is equivalent to the following statement.
For any triangulation $T\in\mathcal{T}_{n}$, we have
\begin{equation}
\bigcup_{\varepsilon\in\{-,+\}^{n}}\left|  \Sigma(T_{\varepsilon})\right|
=\mathcal{T}_{n}. \label{eg}%
\end{equation}
In other words, for any pair $T,T^{\prime}$ of triangulations in
$\mathcal{T}_{n}$, there exist $\varepsilon, \varepsilon^{\prime}%
\in\{-,+\}^{n}$ and a sequence of signed flips from $T_{\varepsilon}$ to
$T^{\prime}_{\varepsilon^{\prime}}$.
\end{theorem}

\subsection{Signed permutations and the four color theorem}

Representing triangulations by permutations via the map $\varphi: S_{n}
\rightarrow\mathcal{T}_{n}$, and using the preceding theorem, we obtain in
this section a reformulation of the four color theorem using now signed permutations.

Consider the alphabet $I=\{\overline{n},\ldots,\overline{1},1,\ldots,n\}$ and
set $I_{-}=\{\overline{n},\ldots,\overline{1}\},$ $I_{+}=\{1,\ldots,n\}.\;$The
letters of $I_{-}$ (resp.\ $I_{+}$) are said \textit{negative}
(resp.\ \textit{positive}).\ Define the bar involution on the letters $\beta$
of $I_{n}$ by
\[
\overline{\beta}=k\text{ if }\beta=\overline{k}\in I_{-}\text{ and }%
\overline{\beta}=\overline{k}\text{ if }\beta=k\in I_{+}\text{.}%
\]
Write $\left|  \beta\right|  =\beta$ if $\beta\in I_{+}$ and $\left|
\beta\right|  =\overline{\beta}$ if $\beta\in I_{-}.\;$For any word
$w=\beta_{1}\cdot\cdot\cdot\beta_{n}$ on $I_{n},$ set $\left|  w\right|
=\left|  \beta_{1}\right|  \cdot\cdot\cdot\left|  \beta_{n}\right|  .\;$The
set of signed permutations on $I_{n}$ is defined by
\[
SP(n)=\{w=\beta_{1}\cdot\cdot\cdot\beta_{n}\mid\left|  \beta_{i}\right|
\neq\left|  \beta_{j}\right|  \text{ for any }i\neq j\}.
\]
Two words $w_{1}$ and $w_{2}$ of $SP(n)$ differ by an authorized transposition
if one of the two situations happens:

\begin{enumerate}
\item $\left|  w_{1}\right|  $ and $\left|  w_{2}\right|  $ are sylvester
adjacent, that is there exist letters $\alpha,\beta,\gamma$ such that $\left|
\alpha\right|  <\left|  \beta\right|  <\left|  \gamma\right|  $ and
%
%
%
%
\begin{equation}%
\begin{array}
[c]{rcl}%
w_{1} & = & u\ \alpha\ \gamma\ v\ \beta\ w\\
w_{2} & = & u\ \gamma\ \alpha\ v\ \beta\ w \label{k1}%
\end{array}
\end{equation}
where $u,v,w$ are factors of $w_{1}$.

\item  There exist letters $\alpha,\gamma$ with the same sign such that
%
%
%
%
\begin{equation}%
\begin{array}
[c]{rcl}%
w_{1} & = & u\ \alpha\ \gamma\ v\\
w_{2} & = & u\ \overline{\gamma}\ \overline{\alpha}\ v\label{k2}%
\end{array}
\end{equation}
where $u,v$ are factors of $w_{1}$ and $v$ does not contain any letter $\beta$
such that $\left|  \alpha\right|  <\left|  \beta\right|  <\left|
\gamma\right|  $ (i.e. $\left|  w_{1}\right|  $ and $\left|  w_{2}\right|  $
are not sylvester adjacent).
\end{enumerate}

Consider $\sigma$ and $\sigma^{\prime}$ two permutations of $S_{n}.\;$We will
say that there exists a signed path between $\sigma$ and $\sigma^{\prime}$ if
one can find two words $w,w^{\prime}$ in $SP(n)$ such that $\left|  w\right|
=\sigma,$ $\left|  w^{\prime}\right|  =\sigma^{\prime}$ and $w^{\prime}$ can
be obtained from $w$ by applying successive authorized transpositions. From
Proposition \ref{prop_sylvSn} and Theorem \ref{ref1} we deduce the following
reformulation of the four color theorem:

\begin{theorem}
\label{fc}The four color theorem is equivalent to the following statement:

\noindent For any positive integer $n$, there exists at least a signable path
joining two permutations of $S_{n}$.
\end{theorem}

\begin{proof}
To each signed triangulation $T_{\varepsilon}$ with $\varepsilon\in
\{-,+\}^{n},$ we associate the set of signed permutations
\[
\lbrack T_{\varepsilon}]_{s}=\{w=\beta_{1}\cdot\cdot\cdot\beta_{n}\mid\left|
w\right|  \in\lbrack T]\text{ and for each }\diamondsuit\in\{-,+\},\beta
_{i}\in I_{\diamondsuit}\Longleftrightarrow\varepsilon_{i}=\diamondsuit\}
\]
that is, $[T_{\varepsilon}]_{s}$ is the set of signed permutations obtained by
signing the words contained in the sylvester class of $T$ so that the sign
associated to each $x_{i}\in X$ is equal to $\varepsilon_{i}.$ By Proposition
\ref{prop_sylvSn}, two signed permutations $w^{\prime}$ and $w$ belong to
$[T_{\varepsilon}]_{s}$ if and only if they differ by successive
transpositions of kind (\ref{k1}).

Given two triangulations $T$ and $T^{\prime}$ in $\mathcal{T}_{n},$ we said
that there exists a signable flip path between $T$ and $T^{\prime}$ if
$T^{\prime}$ belongs to $\bigcup_{\varepsilon\in\{-,+\}^{n}}\left|
\Sigma(T_{\varepsilon})\right|  $ (with the notation of \ref{subsec_fourcol}).
We deduce from Proposition \ref{prop_flip} that the signed triangulations
$T_{\varepsilon}$ and $T_{\varepsilon^{\prime}}^{\prime}$ differ by a signed
flip if and only if there exist a signed word in $[T_{\varepsilon}]_{s}$ and a
signed word in $[T_{\varepsilon^{\prime}}^{\prime}]_{s}$ which differ by a
transposition of kind (\ref{k2}).\ Thus, $T_{\varepsilon}$ and $T_{\varepsilon
^{\prime}}^{\prime}$ differ by a signed flip if and only if each signed word
of $[T_{\varepsilon^{\prime}}^{\prime}]_{s}$ can be obtained by applying
successive transpositions of kind (\ref{k1}) or (\ref{k2}) from any signed
word of $[T_{\varepsilon}]_{s}.$ By transitivity, we deduce that the existence
of a signable flip path between $T$ and $T^{\prime}$ is equivalent to that of
a signable path between each reading of $T$ and each reading of $T^{\prime}.$

Now, by Theorem \ref{ref1} the four color theorem is equivalent to the statement:

\noindent For any positive integer $n$, there exists at least a signable flip
path joining two triangulations of $\mathcal{T}_{n}$.By the previous argument,
this implies our theorem.
\end{proof}

\bigskip

\noindent\textbf{Example.} Consider
\[%
\begin{array}
[c]{rcl}%
\sigma_{1} & = & 324156,\\
\sigma_{2} & = & 453126.
\end{array}
\]
Here is a signed path between suitable signings of $\sigma_{1}$ and
$\sigma_{2}$, where $\rightarrow$ indicates a signed flip and $\equiv$ a
sylvester adjacency:
\[
\bar{3}2\bar{4}156 \rightarrow\bar{3}2\bar{4}\bar{5}\bar{1}6 \rightarrow
\bar{3}254\bar{1}6 \equiv\bar{3}524\bar{1}6 \rightarrow\bar{3}5\bar{4}\bar
{2}\bar{1}6 \equiv5\bar{3}\bar{4}\bar{2}\bar{1}6 \rightarrow543\bar{2}\bar{1}6
\rightarrow\bar{4}\bar{5}3\bar{2}\bar{1}6 \rightarrow\bar{4}\bar{5}3126.
\]
This signed path produces an explicit Heawood signing, and hence an explicit
proper four vertex-coloring, of the sphere triangulation obtained by gluing
the octogon triangulations $\varphi(324156)$ and $\varphi(453126)$ along their boundary.

\bigskip

\noindent\textbf{Remark.} There exists a simple procedure deciding whether a
given path $(\sigma_{1},\ldots,\sigma_{r})$ in $\mathrm{Cay}(S_{n})$ is
signable. The path $(\sigma_{1},\ldots,\sigma_{r})$ is signable if one can
compute a sequence $w_{2},\ldots,w_{r}$ of signed permutations such that
$\left|  w_{i}\right|  =\sigma_{i}$ for any $i=2,\ldots,r$ by the following
procedure. First, the permutations $\sigma_{1}$ and $\sigma_{2}$ differ by a
transposition, thus can be written
\[%
\begin{array}
[c]{rcl}%
\sigma_{1} & = & u\ a\ c\ v\\
\sigma_{2} & = & u\ c\ a\ v.
\end{array}
\]
where $a,c\in\{1,\ldots,n\}$.\ If $\sigma_{1}$ and $\sigma_{2}$ differ by a
sylvester relation, set $w_{2}=\sigma_{2}.\;$Otherwise set $w_{2}%
=u\ \overline{c}\ \overline{a}\ v$.\ Now suppose we have obtained the signed
permutations $w_{2},\ldots,w_{i}$ from $\sigma_{2},\ldots,\sigma_{i}$.\ Since
$\sigma_{i}$ and $\sigma_{i+1}$ differ by a transposition we have%
\[%
\begin{array}
[c]{lcl}%
w_{i} & = & U_{i}\ \alpha_{i}\ \gamma_{i}\ V_{i}\\
\sigma_{i+1} & = & u_{i}\ c_{i}\ a_{i}\ v_{i}%
\end{array}
\]
where $\left|  \alpha_{i}\right|  =a_{i},\left|  \gamma_{i}\right|  =c_{i}$
and $\left|  U_{i}\right|  =u_{i},\left|  V_{i}\right|  =v_{i}$.\ If
$\sigma_{i}$ and $\sigma_{i+1}$ differ by a sylvester relation, set
$w_{i+1}=w_{i}.\;$Otherwise two situations can happen.

\begin{itemize}
\item  When $\alpha_{i}$ and $\gamma_{i}$ have opposite signs the algorithm
stops and the path is not signable.

\item  When $\alpha_{i}$ and $\gamma_{i}$ have the same sign, set
$w_{i+1}=U_{i}\ \overline{\gamma}_{i}\ \overline{\alpha}_{i}\ V_{i}$.
\end{itemize}

\noindent\textbf{Further remarks.}

\noindent$\mathrm{(i)}$ To decide if there exists a signed path between the
two permutations $\sigma$ and $\sigma^{\prime},$ it is sufficient by Lemma
\ref{lem-esse} to test the paths joining these vertices in $\mathrm{Cay}%
(S_{n})$ by restricting to the paths with no loop.

\noindent$\mathrm{(ii)}$ If $(\sigma_{1},\ldots,\sigma_{r})$ is signable in
$\mathrm{Cay}(S_{n}),$ then $(T_{1},\ldots,T_{r})$ where for any
$i=1,\ldots,r,$ $T_{i}=[\sigma_{i}]$ is a signable path in $\mathcal{F}_{n}$.

\subsection{Diagonal signings}

We are going to describe an alternative formulation of this theorem using
signings of the \textit{diagonals} in the triangulations, rather than of the
faces, to be used later in subsection \ref{subsec_asso}.\ Consider a signed
triangulation $T_{\varepsilon}$ and $d$ a diagonal of $T.\;$Let $\varepsilon
_{i}$ and $\varepsilon_{j}$ be the signs of the two faces of $T_{\varepsilon}$
adjacent to $d.$ Then label the diagonal $d$ by the product $\varepsilon
_{i}\varepsilon_{j}\in\{-,+\}.\;$By proceeding similarly for all the diagonals
of $T,$ we obtain a triangulation with signed diagonals $T_{\delta}$. Note
that the signing $T_{\delta}$ of the diagonals of $T$ determines the signing
of its faces up to an overall inversion of signs.\ Indeed, the signing of the
faces of $T$ can be obtained from the sign of one face and the signs of the
diagonals of $T$.

Suppose that the diagonal $d$ in $T$ is flipped in $d^{\prime}$ and write $Q$
for the quadrilateral whose diagonals are $d$ and $d^{\prime}$. The signed
flip of the diagonal $d$ is defined in $T_{\delta}$ when the sign attached to
$d$ is $+.\;$In this case, this yields a triangulation with signed diagonals
defined by the three followings requirements:

\begin{enumerate}
\item $d^{\prime}$ is signed by $+$

\item  the signs of the diagonals of $T_{\delta}$ which are edges of $Q$ are
changed into their opposite

\item  the signs of the remaining diagonals are unchanged.
\end{enumerate}

Note that the signed flip is not defined when the sign of $d$ is $-$. Denote
by $\Sigma(T_{\delta})$ the subgraph of the flip graph generated from
$T_{\delta}$ by applying signed flips and write $\left|  \Sigma(T_{\delta
})\right|  $ for the set of triangulations obtained by deleting the signs $-$
and $+$ in the triangulations with signed diagonals belonging to
$\Sigma(T_{\delta}).\;$By definition of the signed flips on triangulations
with signed diagonals, we have $\left|  \Sigma(T_{\varepsilon})\right|
=\left|  \Sigma(T_{\delta})\right|  $ and Theorem \ref{ref1} is equivalent to
the assertion:

\noindent For any triangulation $T\in\mathcal{T}_{n}$%
\[
\bigcup_{\delta\in\{-,+\}^{n-1}}\left|  \Sigma(T_{\delta})\right|
=\mathcal{T}_{n}%
\]
where $\delta$ yields a signing $T_{\delta}$ of the $n-1$ diagonals of $T$.

\bigskip

We suppose in the sequel that we have chosen a labelling of the $\frac
{n(n+1)}{2}$ diagonals of $P$ (for example we can consider the labelling by
roots belonging to the root system of type $A_{n-1}$ depicted in
\cite{FR}).\ Consider a path $\mathcal{P=}(T_{1},\ldots,T_{r})$ in
$\mathcal{F}_{n}$ where for $i=1,\ldots,r-1,$ the triangulations $T_{i}$ and
$T_{i+1}$ differs by a diagonal flip. The path $\mathcal{P}$ is said signable
if for each $i\in\{1,\ldots,r\},$ there exists a signing $T_{\delta^{(i)}}$ of
the diagonals in $T_{i}$ such that the transformation $T_{\delta^{(i)}%
}\rightarrow T_{\delta^{(i+1)}}$ is a signed flip (that is, Conditions $1,2$
and $3$ above are verified). In fact it is rather easy to determinate if a
path is signable or not.\ For any $i=2,\ldots,r,$ write $N_{i}$ for the set of
diagonals in $T_{i}$ which has appeared after a flip $T_{k}\rightarrow
T_{k+1},$ $1\leq k\leq i-1$. This means that $N_{i}$ is the set of diagonals
of $T_{i}$ which result of a flip operation at one step of the path
$T_{1},\ldots,T_{i}$.\ 

Consider first the flip $T_{1}\rightarrow T_{2}$ of the diagonal $d$ into
$d^{\prime}$. We have $N_{2}=\{d^{\prime}\}$ and the sign of $d^{\prime}$ in
$T_{2}$ must be $+.\;$Now suppose by induction that we have determined the
signs of the diagonals of $N_{i}$ $i\geq2$ so that at each step $T_{k}%
\rightarrow T_{k+1},$ $1\leq k\leq i-1,$ the signs obtained are compatible
with a signed flip (that is with Conditions $1,2$ and $3$ above).\ Let $d_{i}$
be the diagonal flipped in $T_{i}$ and $d_{i}^{\prime}$ the new diagonal
obtained in $T_{i+1}.\;$When $d_{i}$ has sign $-,$ the path is not
signable.\ When $d_{i}$ is not signed or is signed by $+,$ we have
$N_{i+1}=N_{i}+\{d_{i}^{\prime}\}-\{d_{i}\}$.\ In $T_{i+1},$ we sign
$d_{i}^{\prime}$ with $+$ and we change the signs of the edges of the
quadrilateral associated to $d_{i}$ and $d_{i}^{\prime}$ which belong to
$N_{i}.$

The path $\mathcal{P}$ is signable if the previous algorithm does not stop
until the last flip $T_{r-1}\rightarrow T_{r}$ has been considered. In this
case, it becomes immediate to complete in each $T_{i}$ the signs of the
diagonals which are not in $N_{i}.$ Indeed, its suffices to sign arbitrary the
unsigned diagonals of $T_{r}$.$\;$The complete signings in $T_{1}%
,\ldots,T_{r-1}$ are then determined by conditions $1,2$ and $3$ above.%

\begin{center}
\includegraphics[
height=6.7063cm,
width=9.8013cm
]%
{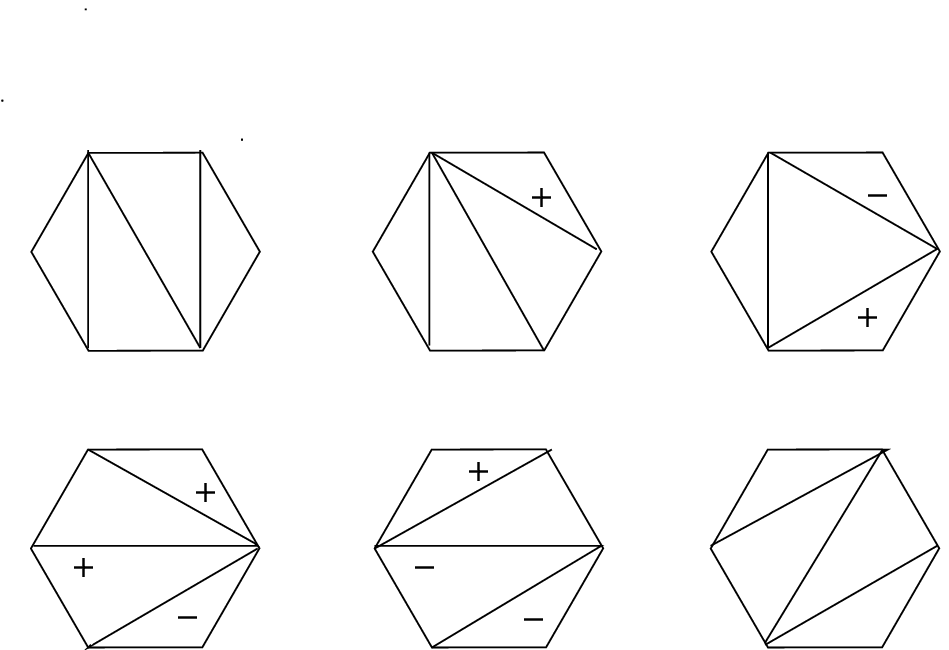}%
\\
Figure 7: Signable and nonsignable paths in $\mathcal{F}_{n}$
\label{fig11}%
\end{center}

Given a path $\mathcal{P=}(T_{1},\ldots,T_{r}),$ we denote by $D_{\mathcal{P}%
}$ the subset of diagonals on $P$ which appear in the triangulations
$T_{1},\ldots,T_{r}$.\ We endow $D_{\mathcal{P}}$ with the structure of an
oriented graph by drawing an edge $d\rightarrow d^{\prime}$ between the
diagonals $d$ and $d^{\prime}$ if there exists $i\in\{1,\ldots,r-1\}$ such
that $T_{i+1}$ is obtained by flipping $d\in T_{i}$ into $d^{\prime}$.

\bigskip

\noindent\textbf{Remarks.}

\noindent$\mathrm{(i)}$ If $\mathcal{P=}(T_{1},\ldots,T_{r})$ is a signable
path, the paths $\mathcal{P=}(T_{p},\ldots,T_{q}),$ $2\leq p<q\leq r$ are also signable.

\noindent$\mathrm{(ii)}$ If all the chains in $D_{\mathcal{P}}$ have length
less than or equal to $2$ (i.e. each diagonal is flipped at most on time), the
path $\mathcal{P=}(T_{1},\ldots,T_{r})$ is signable.\ Indeed $r=1$ or for any
$i\in\{2,\ldots,r\},$ the diagonal $d_{i}$ in the above procedure is never signed.

\noindent$\mathrm{(iii)}$ Suppose that $T$ and $U$ are two triangulations of
$\mathcal{F}_{n}$.\ By Lemma \ref{lem-esse} the problem of determining whether
there exists a signable diagonal path between $T$ and $U$ can be solved by
considering only paths without loop between $T$ and $U$ in $\mathcal{F}_{n}%
$.\ There is a finite number of such paths and one can apply to each of them
the previous procedure.

\bigskip

>From the above arguments and Theorem \ref{ref1} we obtain the following
reformulation of the four color theorem:

\begin{corollary}
\label{cor}(of Theorem \ref{ref1})

\noindent The four color theorem is equivalent to the following statement:

\noindent For any positive integer $n,$ there exists at least a signable path
between two triangulations of the $(n+2)$-gon.
\end{corollary}

\subsection{Signed walks on the associahedron\label{subsec_asso}}

The flip graph $\mathcal{F}_{n}$ is the $1$-skeleton of a convex polytope
called the $n$-dimensional associahedron $\frak{A}_{n}$ (see \cite{FR}). The
vertices of the associahedron can be identified with the triangulations of the
$(n+2)$-gon $P$ and its facets with the $\frac{n(n+1)}{2}$ diagonals of
$P.\;$Its edges correspond to partial triangulations of $P$.\ Given a star
$S,$ write $v_{S}$ for the vertex of maximal order in $S$. This yields a
natural one-to-one correspondence between the stars and the vertices of
$P.\;$The vertices of $\frak{A}_{n}$ belonging to the facet $F_{d}$
corresponding to the diagonal $d$ coincide with the triangulations which
contains $d$.\ This implies that each facet $F_{d}$ contains exactly two
stars. More precisely we have:

\begin{lemma}
\label{lem_facet}Let $S_{1}$ and $S_{2}$ be two stars considered as vertices
of $\frak{A}_{n}$.\ Then we have the following equivalences:

\begin{enumerate}
\item $S_{1}$ and $S_{2}$ belong to the same facet $F_{d}$ if and only if
$v_{S_{1}}v_{S_{2}}=d.$

\item $S_{1}$ and $S_{2}$ do not belong to the same facet $F_{d}$ if and only
if $v_{S_{1}}v_{S_{2}}$ are two consecutive vertices of $P$ (that is form an
edge of $P$).
\end{enumerate}
\end{lemma}

Consider a triangulation $T$ and a diagonal $d$ in $T$. Denote by $T^{\prime}$
the triangulation obtained by flipping $d$ in $T$ and by $d^{\prime}$ the
diagonal of $T^{\prime}$ such that $T^{\prime}=T-\{d\}+\{d^{\prime}\}$.\ The
diagonal flip $T\rightarrow T^{\prime}$ can be interpreted as a move from
$F_{d}$ to $F_{d^{\prime}}$ in $\frak{A}_{n}$ along the edge defined by the
partial triangulation $T-\{d,d^{\prime}\}$.$\frak{\ }$In this case, the facets
$F_{d}$ and $F_{d^{\prime}}$ are disjoint. Indeed the diagonals $d$ and
$d^{\prime}$ are secant in $P,$ thus one cannot find a triangulation $T$
belonging to $F_{d}$ and $F_{d^{\prime}}$.

An edge $E$ of $\frak{A}_{n}$ is contained in $n-2$ facets.\ Suppose that
$E=TT^{\prime}$ where $T$ and $T^{\prime}$ are triangulations belonging
respectively to $F_{d}$ and $F_{d^{\prime}}.\;$Then $T^{\prime}$ is obtained
from $T$ by flipping $d$ into $d^{\prime}.\;$Denote by $Q$ the quadrilateral
in $P$ whose diagonals are $d$ and $d^{\prime}.\;$The edges of $Q$ are either
diagonals either edges of $P$.\ Moreover, for $n\geq2,$ the number of edges of
$Q$ which are diagonals of $P$ belongs to $\{1,2,3,4\}.\;$To distinguish the
facets corresponding to diagonals of $P$ among the facets containing $E,$ it
suffices to consider the four stars $S_{1},S_{2}$ and $S_{3},S_{4}$ which
appear respectively in $F_{d}$ and $F_{d^{\prime}}.\;$Then the facets
corresponding to diagonals of $P$ are the diagonals which can be obtained by
connecting two vertices $S_{i}S_{j},$ $i\neq j.\;$By Lemma \ref{lem_facet},
this is equivalent to find the pairs of vertices $(S_{i},S_{j})$ in
$\frak{A}_{n}$ which belongs to the same facet. The pairs $(S_{1},S_{2})$ and
$(S_{3},S_{4})$ correspond respectively to $d$ and $d^{\prime}$.\ In the
sequel we will denote by $\Delta_{E}$ the set of facets corresponding to the
edges of $Q$ which are diagonals of $P$ distinct of $d$ and $d^{\prime}$. By
the previous arguments $\Delta_{E}$ coincide with the facets containing a pair
$(S_{i},S_{j})$ distinct of $(S_{1},S_{2})$ and $(S_{3},S_{4})$ where
$S_{1},S_{2},S_{3},S_{4}$ are the stars appearing in $F_{d}$ and
$F_{d^{\prime}}$ the two facets of $\frak{A}_{n}$ connected by $E$.

\bigskip

Our aim is now to obtain a reformulation of the four color theorem using only
the geometry of the associahedron and the distinguish subset $\frak{S}_{n}$ of
its vertices which correspond to stars. A walk $W$ of length $r$ on
$\frak{A}_{n}$ is determined by $r$ vertices $V_{1},\ldots,V_{r}$ such that
for any $i=1,\ldots,r-1,$ $V_{i}V_{i+1}$ is an edge of $\frak{A}_{n}.\;$In
this case write $W=(V_{1},\ldots,V_{r}).$ For each edge $V_{i}V_{i+1}$, denote
by $(F_{i},F_{i}^{\prime}$ $)$ the unique pair of facets in $\frak{A}_{n}$
such that $V_{i}\in F_{i},$ $V_{i+1}\in F_{i}^{\prime}$ and $F_{i}\cap
F_{i}^{\prime}=\emptyset.\;$Set $M_{i+1}=\{F_{1}^{\prime},\ldots,F_{i}%
^{\prime}\},$ $i=1,\ldots,r-1.$ Note that $F_{i}^{\prime}\neq F_{i+1}$ in
general (see example below).

The walk $W$ is signable if the facets of $M_{r}$ can be signed by the
following recursive procedure. First sign the face $F_{1}^{\prime}$ with
$+.\;$Suppose that the facets of $M_{i}$ are signed.\ If $F_{i}\in M_{i}$ is
signed by a $-,$ the algorithm stops.\ Otherwise, consider the edge
$V_{i}V_{i+1}$ and the set of facets $\Delta_{i+1}=\Delta_{V_{i}V_{i+1}}$.
Then sign the face $F_{i}^{\prime}$ with a $+$ and change the sign of the
facets of $\Delta_{i+1}\cap M_{i}.$

>From Corollary \ref{cor} and the above arguments, we derive the following
reformulation of the four color theorem in terms of the geometry of the associahedron:

\begin{theorem}
The four color theorem is equivalent to the following statement:

\noindent For any positive integer $n,$ there is a least a signable walk
between two vertices of the associahedron $\frak{A}_{n}$.
\end{theorem}

\noindent\textbf{Remarks.}

\noindent$\mathrm{(i)}$ By Lemma \ref{lem-esse}, the problem of finding a
signed walk between two vertices $V$ and $V^{\prime}$ of the associahedron can
be solved by applying the previous procedure to the walks joining $V$ to
$V^{\prime}$ in which each vertex is attained at more one time (that is by
excluding the walks with loops).

\noindent$\mathrm{(ii)}$ At each step of the above procedure, the set
$\Delta_{i+1}=\Delta_{V_{i}V_{i+1}}$ is determined only by the vertices of
$F_{i}$ and $F_{i+1}$ which belong to $\frak{S}_{n}$.

\begin{example}
The signable path of Figure 7 is equivalent to the signed walk on
$\frak{A}_{4}$ given in the figure below%
\[%
\raisebox{-0cm}{\parbox[b]{10.2692cm}{\begin{center}
\includegraphics[
height=12.1166cm,
width=10.2692cm
]%
{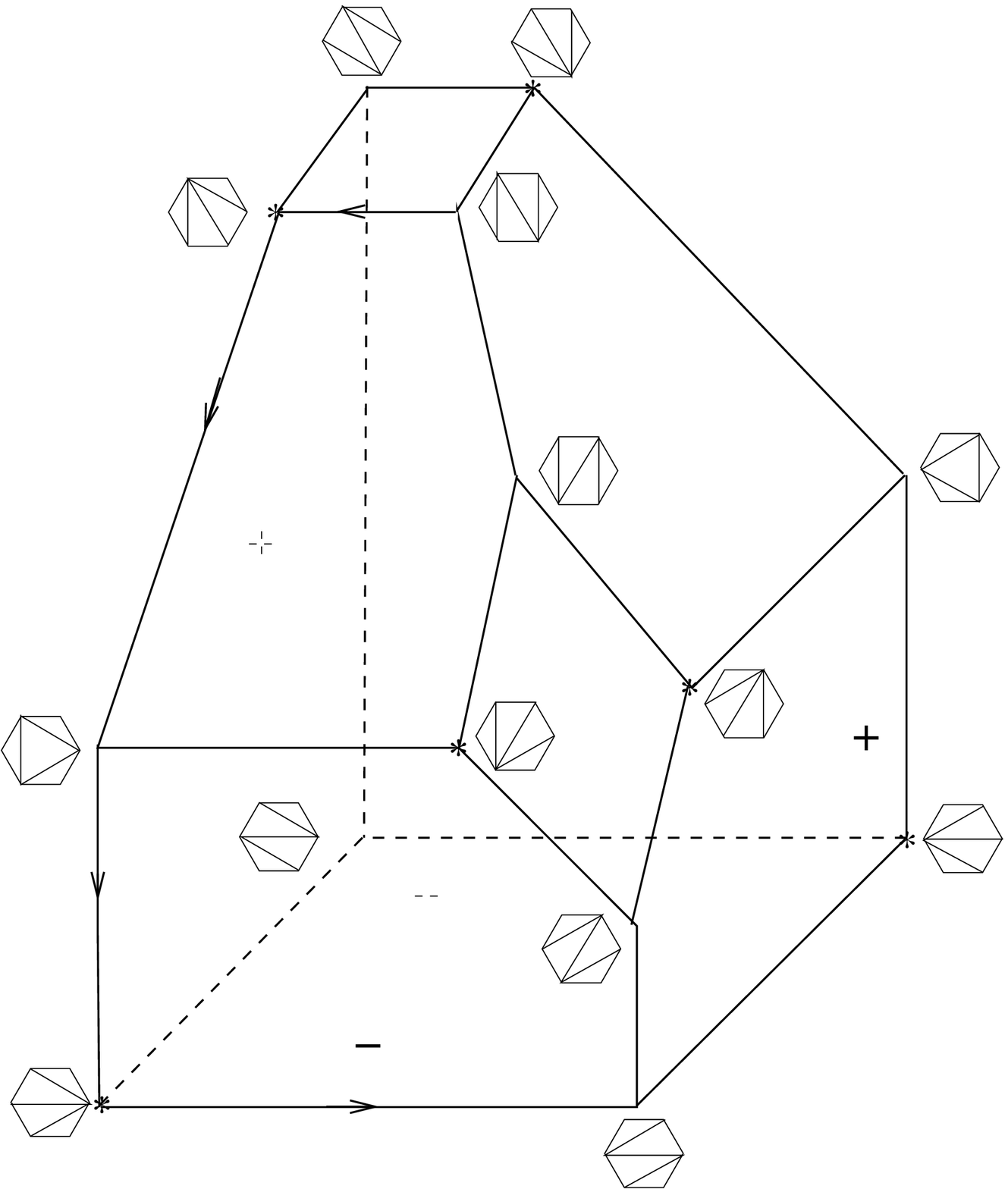}%
\\
Figure 8
\end{center}}}%
\]
where the solid (resp. dashed) signs belong to apparent (resp. non apparent) facets.
\end{example}

\section{Flip graph generated by a colored triangulation \label{sec_fourcolor}}

\subsection{A combinatorial problem}

The set $\mathcal{T}_{n}$ can be identified with the set of simple colored
triangulations $\mathcal{ST}_{n}(\varepsilon)$ with $\varepsilon
=(1,\ldots,1)\in\mathbb{N}^{n}$.\ The flip graph can be generated starting
from any triangulation, by applying diagonal flips.\ This yields to the
following natural problem:

\begin{problem}
\label{problem}What is the graph generated from a colored triangulation by
applying successive restrictive flips?
\end{problem}

\subsection{Homogeneous flips case}

Consider a colored triangulation $T_{\varepsilon}$ and denote by
$\mathcal{HF}(T_{\varepsilon})$ the subgraph of the flip graph generated from
$T_{\varepsilon}$ by applying homogeneous flips. The coloring $\varepsilon$
defines subtriangulations in $T_{\varepsilon}$ obtained by gluing together the
adjacent faces which have the same color. These subtriangulations will be
called the connected components of $T_{\varepsilon}.\;$Denote them by
$T_{1},\ldots,T_{r}$.\ For any $i=1,\ldots,r,$ let $\nu_{i}$ be the number of
faces in $T_{i}.$ Since the homogeneous flips stabilize the connected
components $T_{i}$, we obtain:

\begin{proposition}
The subgraph $\mathcal{HF}(T_{\varepsilon})$ is isomorphic to the direct
product of flip graphs $\mathcal{F}_{\nu_{1}}\times\cdot\cdot\cdot
\times\mathcal{F}_{\nu_{r}}.$
\end{proposition}

The above isomorphism can be explicited by associating to each colored
triangulation of $\mathcal{HF}(T_{\varepsilon})$ the $r$-uple of
triangulations defined from its connected components $T_{1},\ldots,T_{r}$ as
pictured in the figure below.\ This answers Problem \ref{problem}.%

\begin{center}
\includegraphics[
height=4.1684in,
width=4.7089in
]%
{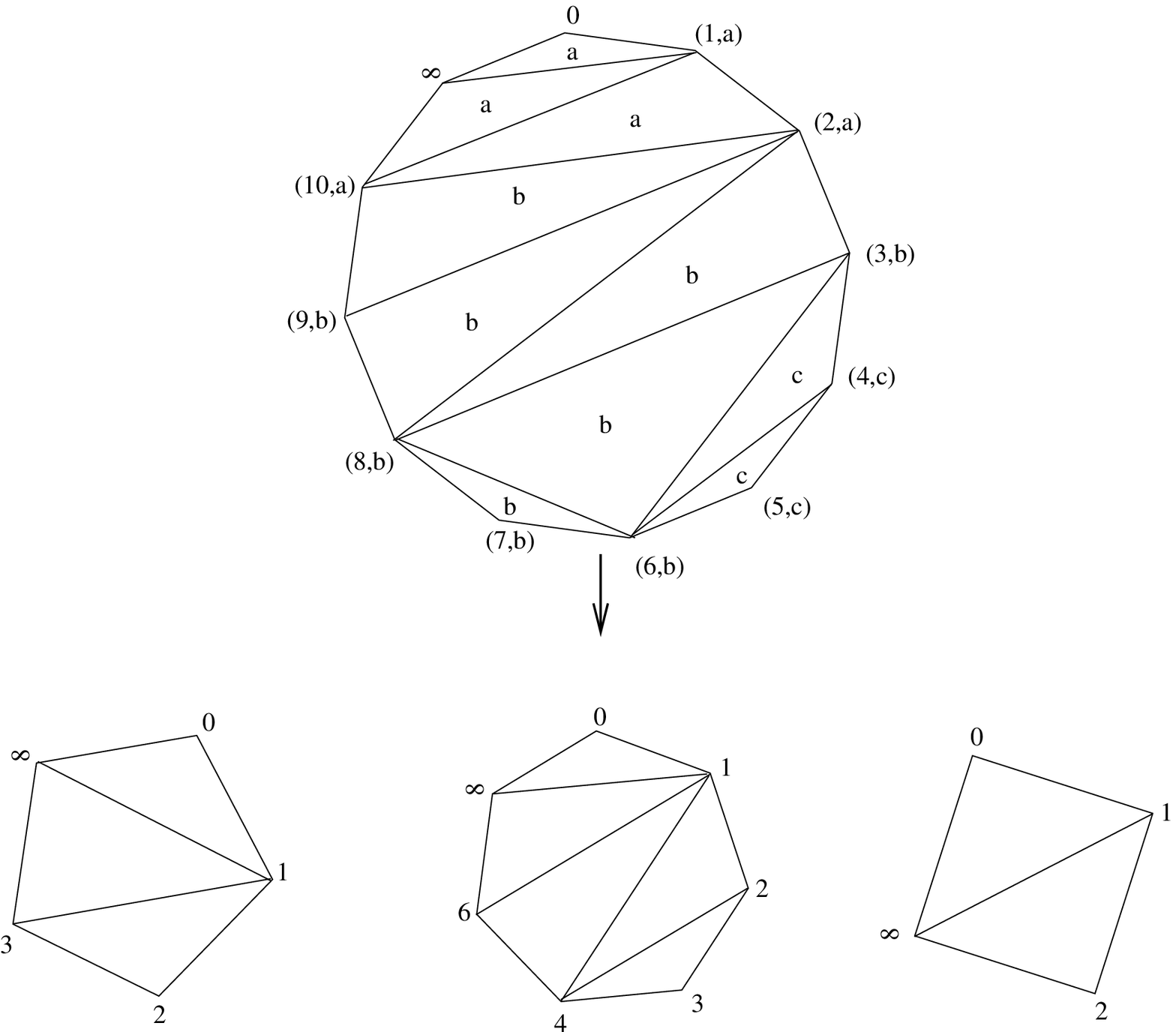}%
\\
Figure 9
\label{fig7}%
\end{center}

\subsection{Switched flips case}

Consider $\mu\in\mathbb{N}^{d}$ and $T_{\varepsilon_{\mu}}$ a simple colored
triangulation of $\mathcal{T}_{n}(\varepsilon_{\mu})$.$\;$Denote by
$\mathcal{SF}(T_{\varepsilon})$ the subgraph of the flip graph generated from
$T_{\varepsilon}$ by applying switched flips. Write $\mathcal{SF}%
_{n}(\varepsilon_{\mu})$ for the graph obtained by drawing an edge between two
simple triangulations of $\mathcal{F}_{n}(\varepsilon_{\mu})$ when they differ
by a switched flip. We are going to show that $\mathcal{SF}_{n}(\varepsilon
_{\mu})$ is connected this will imply the equality $\mathcal{SF}%
(\varepsilon_{\mu})=\mathcal{SF}(T_{\varepsilon})$ for any simple colored
triangulation such that $\varepsilon=\varepsilon_{\mu}$.

Let $S_{\mu}$ be the Frobenius subgroup of $S_{n}$ defined by $\mu,$ that is
the subgroup of permutations which stabilize the intervals $I_{\mu_{1}%
}=\{1,\ldots,\mu_{1}\}$ and $I_{\mu_{s}}=\{\mu_{s-1}+1,\ldots,\mu_{s}\}$ for
$s=2,\ldots,d.$ The elements of the coset $S_{n}/S_{\mu}$ will be identified
with the words of length $n$ and evaluation $\mu$ on the totally ordered
alphabet $\mathcal{C}=\{c_{1}<\cdot\cdot\cdot<c_{d}\},$ that is we set
$S_{n}/S_{\mu}=\mathcal{C}_{n,\mu}$.

Denote by $\mathrm{Cay}(\mathcal{C}_{n,\mu})$ the Cayley graph of
$\mathcal{C}_{n,\mu}.$ This means that the vertices of $\mathcal{C}_{n,\mu}$
are the words of length $n$ and evaluation $\mu$ and there is an edge between
$w$ and $w^{\prime}$ if and only if $w^{\prime}$ is obtained by switching two
adjacent letters of $w.$ Denote by $\mathrm{STD}$ the standardization map on
simple colored triangulations defined by $\mathrm{STD(}T_{\varepsilon_{\mu}%
})=T.$ From the definition of the standardization map and since switched flips
are particular cases of flips, we have:

\begin{lemma}
\label{lemstd} \ \ \ 

\begin{enumerate}
\item  The standardization map $\mathrm{std}$ on words is an injective
morphism of graphs from $\mathrm{Cay}(\mathcal{C}_{n,\mu})$ to $\mathrm{Cay}%
(S_{n})$.

\item  The standardization map $\mathrm{STD}$ is an injective morphism of
graphs from $\mathcal{SF}_{n}(\varepsilon_{\mu})$ to $\mathcal{F}_{n}$.
\end{enumerate}
\end{lemma}

\begin{proposition}
\label{sflip}Let $T_{\varepsilon_{\mu}}$ and $T_{\varepsilon_{\mu}}^{\prime}$
two triangulations in $\mathcal{SF}_{n}(\varepsilon_{\mu})$.\ Then
$T_{\varepsilon_{\mu}}$ and $T_{\varepsilon_{\mu}}^{\prime}$ differ by a
switched flip if and only if there exist a reading $w$ of $T_{\varepsilon
_{\mu}}$ and a reading $w^{\prime}$ of $T_{\varepsilon_{\mu}}^{\prime}$ of the
form
\begin{equation}%
\begin{array}
[c]{lcl}%
w & = & uxzv\\
w^{\prime} & = & uzxv \label{sw}%
\end{array}
\end{equation}
where $x,z$ are letters of $\mathcal{C}$ and $u,v$ words on $\mathcal{C}$ such
that $v$ does not contain any letter $y$ verifying $x\leq y<z$ or $z\leq y<x$.
\end{proposition}

\begin{proof}
By symmetry we only consider the case $x\leq y<z$.\ Suppose that $w$ and
$w^{\prime}$ are readings of $T_{\varepsilon_{\mu}}$ and $T_{\varepsilon_{\mu
}}^{\prime}$ verifying $w=uxzv$ and $w^{\prime}=uzxv$ as in the theorem. Then
$\mathrm{std}(w)=u_{\ast}x_{\ast}z_{\ast}v_{\ast}$ and $\mathrm{std}%
(w^{\prime})=u_{\ast}z_{\ast}x_{\ast}v_{\ast}$ where $x_{\ast},z_{\ast}$ are
letters of $X$ and $u_{\ast},v_{\ast}$ words on $X$ such that $v_{\ast}$ does
not contain any letter $y$ verifying $x<y<z.$ Thus by applying $1$ of
Proposition \ref{prop_flip}, $T$ and $T^{\prime}$ differs by a diagonal flip
and the faces corresponding to its flip are labelled by $x_{\ast}$ and
$z_{\ast}.$ This implies that $T_{\varepsilon_{\mu}}$ and $T_{\varepsilon
_{\mu}}^{\prime}$ differ by a flip corresponding to faces colored by $x$ and
$z$, hence by a switched flip.

Conversely, suppose that $T_{\varepsilon_{\mu}}$ and $T_{\varepsilon_{\mu}%
}^{\prime}$ differ by a switched flip.\ Then $T$ and $T^{\prime}$ differ by a
flip and by $2$ of Proposition \ref{prop_flip}, we have readings
$\sigma=u_{\ast}x_{\ast}z_{\ast}v_{\ast}$ and $\sigma^{\prime}=u_{\ast}%
z_{\ast}x_{\ast}v_{\ast}$ respectively of $T$ and $T^{\prime}$ where $x_{\ast
},z_{\ast}$ are letters of $X$ and $u_{\ast},v_{\ast}$ words on $X$ such that
$v_{\ast}$ does not contain any letter $y$ verifying $x<y<z.$ Hence by
applying the destandardization procedure (see \ref{subsec_asso}) to $\sigma$
and $\sigma^{\prime}$, one can find readings $w=uxzv$ and $w^{\prime}=uzxv$
respectively of $T_{\varepsilon_{\mu}}$ and $T_{\varepsilon_{\mu}}^{\prime}$
with $x\leq y<z.$
\end{proof}

\begin{theorem}
\ \ \ \ 

\begin{enumerate}
\item  The following diagram commutes:%
\[%
\begin{array}
[c]{ccc}%
\mathrm{Cay}(\mathcal{C}_{n,\mu}) & \overset{\mathrm{std}}{\rightarrow} &
\mathrm{Cay}(S_{n})\\
\Phi\downarrow &  & \downarrow\varphi\\
\mathcal{SF}_{n}(\varepsilon_{\mu}) & \underset{\mathrm{STD}}{\rightarrow} &
\mathcal{F}_{n}%
\end{array}
\]
where $\varphi,\Phi,\mathrm{std}$ and $\mathrm{STD}$ are morphisms of graphs.

\item $\mathcal{SF}_{n}(\varepsilon_{\mu})$ is connected.
\end{enumerate}
\end{theorem}

\begin{proof}
\ \ 

$1:$ We have already seen that $\Phi$ is a surjective map. Suppose that $w$
and $w^{\prime}$ are words in $\mathrm{Cay}(\mathcal{C}_{n,\mu})$ which differ
by the transposition of two consecutive letters. If $w$ and $w^{\prime}$ are
sylvester adjacent, $\mathrm{std}(w)$ and $\mathrm{std}(w^{\prime})$ are also
sylvester adjacent, thus $\varphi(\mathrm{std}(w))=\varphi(\mathrm{std}%
(w^{\prime}))$.\ This implies that $\Phi(w)=\Phi(w^{\prime})$ since $w$ and
$w^{\prime}$ have the same evaluation.\ If $w$ and $w^{\prime}$ are not
sylvester adjacent, they verify (\ref{sw}) and by Proposition \ref{sflip}, we
know that $\Phi(w)$ and $\Phi(w^{\prime})$ differ by a switched flip in
$\mathcal{SF}_{n}(\varepsilon_{\mu})$.\ This proves that $\Phi$ is a morphism
of graphs. Now for any word $w$ in $\mathrm{Cay}(\mathcal{C}_{n,\mu})$,
$\mathrm{STD}\circ\Phi(w)=\varphi\circ\mathrm{std}(w)$ by definition of the
map $\Phi.\;$Thus by Lemma \ref{lemstd}, the above diagram is a commutting
diagram of morphisms.

\noindent$2:$ Since $\mathrm{Cay}(\mathcal{C}_{n,\mu})$ is connected and
$\Phi$ is a surjective morphism of graphs, we obtain immediately that
$\mathcal{SF}_{n}(\varepsilon_{\mu})$ is connected.
\end{proof}

\bigskip

\noindent\textbf{Remark. }The theorem implies in particular that
$\mathcal{SF}(T_{\varepsilon})=\mathcal{SF}_{n}(\varepsilon_{\mu})$ for any
simple colored triangulation $T_{\varepsilon_{\mu}}\in\mathcal{T}%
_{n}(\varepsilon_{\mu})$.\ Thus it answers to Problem \ref{problem} when the
colored triangulation $T_{\varepsilon}$ is simple. When $T_{\varepsilon}$ is
not simple, we have find no algebraic interpretation of the graph
$\mathcal{SF}(T_{\varepsilon}).$

\end{document}